\documentclass[11pt,a4paper]{article}

\usepackage{amsmath,amssymb,latexsym,pstricks,mathrsfs,comment,amsthm,graphicx,tikz,tikz-cd,enumerate,
pgffor,cite,wrapfig,multicol,float}
\usepackage{bibspacing}
\newcommand{\nc}{\newcommand}
\nc{\rnc}{\renewcommand}

\usepackage{geometry}  \rnc{\ss}{\smallskip} \nc{\ms}{\medskip} \nc{\bs}{\bigskip} \nc{\nss}{\vspace{-3mm}}

\makeatletter

\begin{document}

\nc{\Reg}{\operatorname{Reg}}
\nc{\RP}{\operatorname{RP}}
\nc{\TXa}{\T_X^a}
\nc{\TXA}{\T(X,A)}
\nc{\TXal}{\T(X,\al)}
\nc{\RegTXa}{\Reg(\TXa)}
\nc{\RegTXA}{\Reg(\TXA)}
\nc{\RegTXal}{\Reg(\TXal)}
\nc{\PalX}{\P_\al(X)}
\nc{\EAX}{\E_A(X)}
\nc{\Bb}{\overline{B}}
\nc{\bb}{\overline{\be}}
\nc{\bw}{{\bf w}}
\nc{\bz}{{\bf z}}
\nc{\TASA}{\T_A\sm\S_A}
\nc{\Ub}{\overline{U}}
\nc{\Vb}{\overline{V}}
\nc{\eb}{\overline{e}}
\nc{\EXa}{\E_X^a}
\nc{\oijr}{1\leq i<j\leq r}
\nc{\veb}{\overline{\ve}}
\nc{\bbT}{\mathbb T}
\nc{\Surj}{\operatorname{Surj}}
\nc{\Sone}{S^{(1)}}
\nc{\fillbox}[2]{\draw[fill=gray!30](#1,#2)--(#1+1,#2)--(#1+1,#2+1)--(#1,#2+1)--(#1,#2);}
\nc{\raa}{\rangle_a}
\nc{\Ea}{E_a}
\nc{\ep}{\varepsilon} \nc{\ve}{\eta}
\nc{\IXa}{\I_X^a}
\nc{\RegIXa}{\Reg(\IXa)}
\nc{\JXa}{\J_X^a}
\nc{\RegJXa}{\Reg(\JXa)}
\nc{\IXA}{\I(X,A)}
\nc{\IAX}{\I(A,X)}
\nc{\RegIXA}{\Reg(\IXA)}
\nc{\RegIAX}{\Reg(\IAX)}
\nc{\trans}[2]{\left(\begin{smallmatrix} #1 \\ #2 \end{smallmatrix}\right)}
\nc{\bigtrans}[2]{\left(\begin{matrix} #1 \\ #2 \end{matrix}\right)}
\nc{\lmap}[1]{\mapstochar \xrightarrow {\ #1\ }}
\nc{\EaTXa}{E}

\nc{\gL}{\mathrel{\mathscr L}}
\nc{\gR}{\mathrel{\mathscr R}}
\nc{\gH}{\mathrel{\mathscr H}}
\nc{\gJ}{\mathrel{\mathscr J}}
\nc{\gD}{\mathrel{\mathscr D}}
\nc{\gK}{\mathrel{\mathscr K}}
\nc{\gLa}{\mathrel{\mathscr L^a}}
\nc{\gRa}{\mathrel{\mathscr R^a}}
\nc{\gHa}{\mathrel{\mathscr H^a}}
\nc{\gJa}{\mathrel{\mathscr J^a}}
\nc{\gDa}{\mathrel{\mathscr D^a}}
\nc{\gKa}{\mathrel{\mathscr K^a}}
\nc{\gLh}{\mathrel{\widehat{\mathscr L}}}
\nc{\gRh}{\mathrel{\widehat{\mathscr R}}}
\nc{\gHh}{\mathrel{\widehat{\mathscr H}}}
\nc{\gJh}{\mathrel{\widehat{\mathscr J}}}
\nc{\gDh}{\mathrel{\widehat{\mathscr D}}}
\nc{\gKh}{\mathrel{\widehat{\mathscr K}}}
\nc{\Lh}{\widehat{L}}
\nc{\Rh}{\widehat{R}}
\nc{\Hh}{\widehat{H}}
\nc{\Jh}{\widehat{J}}
\nc{\Dh}{\widehat{D}}
\nc{\Kh}{\widehat{K}}
\nc{\gLb}{\mathrel{\widehat{\mathscr L}}}
\nc{\gRb}{\mathrel{\widehat{\mathscr R}}}
\nc{\gHb}{\mathrel{\widehat{\mathscr H}}}
\nc{\gJb}{\mathrel{\widehat{\mathscr J}}}
\nc{\gDb}{\mathrel{\widehat{\mathscr D}}}
\nc{\gKb}{\mathrel{\widehat{\mathscr K}}}
\nc{\Lb}{\widehat{L}}
\nc{\Rb}{\widehat{R}}
\nc{\Hb}{\widehat{H}}
\nc{\Jb}{\widehat{J}}
\nc{\Db}{\overline{D}}
\nc{\Kb}{\widehat{K}}

\hyphenation{mon-oid mon-oids}

\nc{\itemit}[1]{\item[\emph{(#1)}]}
\nc{\E}{\mathcal E}
\nc{\TX}{\T(X)}
\nc{\TXP}{\T(X,\P)}
\nc{\EX}{\E(X)}
\nc{\EXP}{\E(X,\P)}
\nc{\SX}{\S(X)}
\nc{\SXP}{\S(X,\P)}
\nc{\Sing}{\operatorname{Sing}}
\nc{\idrank}{\operatorname{idrank}}
\nc{\SingXP}{\Sing(X,\P)}
\nc{\De}{\Delta}
\nc{\sgp}{\operatorname{sgp}}
\nc{\mon}{\operatorname{mon}}
\nc{\Dn}{\mathcal D_n}
\nc{\Dm}{\mathcal D_m}

\nc{\lline}[1]{\draw(3*#1,0)--(3*#1+2,0);}
\nc{\uline}[1]{\draw(3*#1,5)--(3*#1+2,5);}
\nc{\thickline}[2]{\draw(3*#1,5)--(3*#2,0); \draw(3*#1+2,5)--(3*#2+2,0) ;}
\nc{\thicklabel}[3]{\draw(3*#1+1+3*#2*0.15-3*#1*0.15,4.25)node{{\tiny $#3$}};}

\nc{\slline}[3]{\draw(3*#1+#3,0+#2)--(3*#1+2+#3,0+#2);}
\nc{\suline}[3]{\draw(3*#1+#3,5+#2)--(3*#1+2+#3,5+#2);}
\nc{\sthickline}[4]{\draw(3*#1+#4,5+#3)--(3*#2+#4,0+#3); \draw(3*#1+2+#4,5+#3)--(3*#2+2+#4,0+#3) ;}
\nc{\sthicklabel}[5]{\draw(3*#1+1+3*#2*0.15-3*#1*0.15+#5,4.25+#4)node{{\tiny $#3$}};}

\nc{\stll}[5]{\sthickline{#1}{#2}{#4}{#5} \sthicklabel{#1}{#2}{#3}{#4}{#5}}
\nc{\tll}[3]{\stll{#1}{#2}{#3}00}

\nc{\mfourpic}[9]{
\slline1{#9}0
\slline3{#9}0
\slline4{#9}0
\slline5{#9}0
\suline1{#9}0
\suline3{#9}0
\suline4{#9}0
\suline5{#9}0
\stll1{#1}{#5}{#9}{0}
\stll3{#2}{#6}{#9}{0}
\stll4{#3}{#7}{#9}{0}
\stll5{#4}{#8}{#9}{0}
\draw[dotted](6,0+#9)--(8,0+#9);
\draw[dotted](6,5+#9)--(8,5+#9);
}
\nc{\vdotted}[1]{
\draw[dotted](3*#1,10)--(3*#1,15);
\draw[dotted](3*#1+2,10)--(3*#1+2,15);
}

\nc{\Clab}[2]{
\sthicklabel{#1}{#1}{{}_{\phantom{#1}}C_{#1}}{1.25+5*#2}0
}
\nc{\sClab}[3]{
\sthicklabel{#1}{#1}{{}_{\phantom{#1}}C_{#1}}{1.25+5*#2}{#3}
}
\nc{\Clabl}[3]{
\sthicklabel{#1}{#1}{{}_{\phantom{#3}}C_{#3}}{1.25+5*#2}0
}
\nc{\sClabl}[4]{
\sthicklabel{#1}{#1}{{}_{\phantom{#4}}C_{#4}}{1.25+5*#2}{#3}
}
\nc{\Clabll}[3]{
\sthicklabel{#1}{#1}{C_{#3}}{1.25+5*#2}0
}
\nc{\sClabll}[4]{
\sthicklabel{#1}{#1}{C_{#3}}{1.25+5*#2}{#3}
}

\nc{\mtwopic}[6]{
\slline1{#6*5}{#5}
\slline2{#6*5}{#5}
\suline1{#6*5}{#5}
\suline2{#6*5}{#5}
\stll1{#1}{#3}{#6*5}{#5}
\stll2{#2}{#4}{#6*5}{#5}
}
\nc{\mtwopicl}[6]{
\slline1{#6*5}{#5}
\slline2{#6*5}{#5}
\suline1{#6*5}{#5}
\suline2{#6*5}{#5}
\stll1{#1}{#3}{#6*5}{#5}
\stll2{#2}{#4}{#6*5}{#5}
\sClabl1{#6}{#5}{i}
\sClabl2{#6}{#5}{j}
}

\nc{\keru}{\operatorname{ker}^\wedge} \nc{\kerl}{\operatorname{ker}_\vee}

\nc{\coker}{\operatorname{coker}}
\nc{\KER}{\ker}
\nc{\N}{\mathbb N}
\nc{\LaBn}{L_\al(\B_n)}
\nc{\RaBn}{R_\al(\B_n)}
\nc{\LaPBn}{L_\al(\PB_n)}
\nc{\RaPBn}{R_\al(\PB_n)}
\nc{\rhorBn}{\rho_r(\B_n)}
\nc{\DrBn}{D_r(\B_n)}
\nc{\DrPn}{D_r(\P_n)}
\nc{\DrPBn}{D_r(\PB_n)}
\nc{\DrKn}{D_r(\K_n)}
\nc{\alb}{\al_{\vee}}
\nc{\beb}{\be^{\wedge}}
\nc{\bnf}{\bn^\flat}
\nc{\Bal}{\operatorname{Bal}}
\nc{\Red}{\operatorname{Red}}
\nc{\Pnxi}{\P_n^\xi}
\nc{\Bnxi}{\B_n^\xi}
\nc{\PBnxi}{\PB_n^\xi}
\nc{\Knxi}{\K_n^\xi}
\nc{\C}{\mathscr C}
\nc{\exi}{e^\xi}
\nc{\Exi}{E^\xi}
\nc{\eximu}{e^\xi_\mu}
\nc{\Eximu}{E^\xi_\mu}
\nc{\REF}{ {\red [Ref?]} }
\nc{\GL}{\operatorname{GL}}
\rnc{\O}{\operatorname{O}}

\nc{\vtx}[2]{\fill (#1,#2)circle(.2);}
\nc{\lvtx}[2]{\fill (#1,0)circle(.2);}
\nc{\uvtx}[2]{\fill (#1,1.5)circle(.2);}

\nc{\Eq}{\mathfrak{Eq}}
\nc{\Gau}{\Ga^\wedge} \nc{\Gal}{\Ga_\vee}
\nc{\Lamu}{\Lam^\wedge} \nc{\Laml}{\Lam_\vee}
\nc{\bX}{{\bf X}}
\nc{\bY}{{\bf Y}}
\nc{\ds}{\displaystyle}

\nc{\uvert}[1]{\fill (#1,1.5)circle(.2);}
\nc{\uuvert}[1]{\fill (#1,3)circle(.2);}
\nc{\uuuvert}[1]{\fill (#1,4.5)circle(.2);}
\rnc{\lvert}[1]{\fill (#1,0)circle(.2);}
\nc{\overt}[1]{\fill (#1,0)circle(.1);}
\nc{\overtl}[3]{\node[vertex] (#3) at (#1,0) {  {\tiny $#2$} };}
\nc{\cv}[2]{\draw(#1,1.5) to [out=270,in=90] (#2,0);}
\nc{\cvs}[2]{\draw(#1,1.5) to [out=270+30,in=90+30] (#2,0);}
\nc{\ucv}[2]{\draw(#1,3) to [out=270,in=90] (#2,1.5);}
\nc{\uucv}[2]{\draw(#1,4.5) to [out=270,in=90] (#2,3);}
\nc{\textpartn}[1]{{\lower0.45 ex\hbox{\begin{tikzpicture}[xscale=.2,yscale=0.2] #1 \end{tikzpicture}}}}
\nc{\textpartnx}[2]{{\lower1.0 ex\hbox{\begin{tikzpicture}[xscale=.3,yscale=0.3] 
\foreach \x in {1,...,#1}
{ \uvert{\x} \lvert{\x} }
#2 \end{tikzpicture}}}}
\nc{\disppartnx}[2]{{\lower1.0 ex\hbox{\begin{tikzpicture}[scale=0.3] 
\foreach \x in {1,...,#1}
{ \uvert{\x} \lvert{\x} }
#2 \end{tikzpicture}}}}
\nc{\disppartnxd}[2]{{\lower2.1 ex\hbox{\begin{tikzpicture}[scale=0.3] 
\foreach \x in {1,...,#1}
{ \uuvert{\x} \uvert{\x} \lvert{\x} }
#2 \end{tikzpicture}}}}
\nc{\disppartnxdn}[2]{{\lower2.1 ex\hbox{\begin{tikzpicture}[scale=0.3] 
\foreach \x in {1,...,#1}
{ \uuvert{\x} \lvert{\x} }
#2 \end{tikzpicture}}}}
\nc{\disppartnxdd}[2]{{\lower3.6 ex\hbox{\begin{tikzpicture}[scale=0.3] 
\foreach \x in {1,...,#1}
{ \uuuvert{\x} \uuvert{\x} \uvert{\x} \lvert{\x} }
#2 \end{tikzpicture}}}}

\nc{\dispgax}[2]{{\lower0.0 ex\hbox{\begin{tikzpicture}[scale=0.3] 
#2
\foreach \x in {1,...,#1}
{\lvert{\x} }
 \end{tikzpicture}}}}
\nc{\textgax}[2]{{\lower0.4 ex\hbox{\begin{tikzpicture}[scale=0.3] 
#2
\foreach \x in {1,...,#1}
{\lvert{\x} }
 \end{tikzpicture}}}}
\nc{\textlinegraph}[2]{{\raise#1 ex\hbox{\begin{tikzpicture}[scale=0.8] 
#2
 \end{tikzpicture}}}}
\nc{\textlinegraphl}[2]{{\raise#1 ex\hbox{\begin{tikzpicture}[scale=0.8] 
\tikzstyle{vertex}=[circle,draw=black, fill=white, inner sep = 0.07cm]
#2
 \end{tikzpicture}}}}
\nc{\displinegraph}[1]{{\lower0.0 ex\hbox{\begin{tikzpicture}[scale=0.6] 
#1
 \end{tikzpicture}}}}
 
\nc{\disppartnthreeone}[1]{{\lower1.0 ex\hbox{\begin{tikzpicture}[scale=0.3] 
\foreach \x in {1,2,3,5,6}
{ \uvert{\x} }
\foreach \x in {1,2,4,5,6}
{ \lvert{\x} }
\draw[dotted] (3.5,1.5)--(4.5,1.5);
\draw[dotted] (2.5,0)--(3.5,0);
#1 \end{tikzpicture}}}}

\nc{\partn}[4]{\left( \begin{array}{c|c} 
#1 \ & \ #3 \ \ \\ \cline{2-2}
#2 \ & \ #4 \ \
\end{array} \!\!\! \right)}
\nc{\partnlong}[6]{\partn{#1}{#2}{#3,\ #4}{#5,\ #6}} 
\nc{\partnsh}[2]{\left( \begin{array}{c} 
#1 \\
#2 
\end{array} \right)}
\nc{\partncodefz}[3]{\partn{#1}{#2}{#3}{\emptyset}}
\nc{\partndefz}[3]{{\partn{#1}{#2}{\emptyset}{#3}}}
\nc{\partnlast}[2]{\left( \begin{array}{c|c}
#1 \ &  \ #2 \\
#1 \ &  \ #2
\end{array} \right)}

\nc{\uuarcx}[3]{\draw(#1,3)arc(180:270:#3) (#1+#3,3-#3)--(#2-#3,3-#3) (#2-#3,3-#3) arc(270:360:#3);}
\nc{\uuarc}[2]{\uuarcx{#1}{#2}{.4}}
\nc{\uuuarcx}[3]{\draw(#1,4.5)arc(180:270:#3) (#1+#3,4.5-#3)--(#2-#3,4.5-#3) (#2-#3,4.5-#3) arc(270:360:#3);}
\nc{\uuuarc}[2]{\uuuarcx{#1}{#2}{.4}}
\nc{\darcx}[3]{\draw(#1,0)arc(180:90:#3) (#1+#3,#3)--(#2-#3,#3) (#2-#3,#3) arc(90:0:#3);}
\nc{\darc}[2]{\darcx{#1}{#2}{.4}}
\nc{\udarcx}[3]{\draw(#1,1.5)arc(180:90:#3) (#1+#3,1.5+#3)--(#2-#3,1.5+#3) (#2-#3,1.5+#3) arc(90:0:#3);}
\nc{\udarc}[2]{\udarcx{#1}{#2}{.4}}
\nc{\uudarcx}[3]{\draw(#1,3)arc(180:90:#3) (#1+#3,3+#3)--(#2-#3,3+#3) (#2-#3,3+#3) arc(90:0:#3);}
\nc{\uudarc}[2]{\uudarcx{#1}{#2}{.4}}
\nc{\uarcx}[3]{\draw(#1,1.5)arc(180:270:#3) (#1+#3,1.5-#3)--(#2-#3,1.5-#3) (#2-#3,1.5-#3) arc(270:360:#3);}
\nc{\uarc}[2]{\uarcx{#1}{#2}{.4}}
\nc{\darcxhalf}[3]{\draw(#1,0)arc(180:90:#3) (#1+#3,#3)--(#2,#3) ;}
\nc{\darchalf}[2]{\darcxhalf{#1}{#2}{.4}}
\nc{\uarcxhalf}[3]{\draw(#1,1.5)arc(180:270:#3) (#1+#3,1.5-#3)--(#2,1.5-#3) ;}
\nc{\uarchalf}[2]{\uarcxhalf{#1}{#2}{.4}}
\nc{\uarcxhalfr}[3]{\draw (#1+#3,1.5-#3)--(#2-#3,1.5-#3) (#2-#3,1.5-#3) arc(270:360:#3);}
\nc{\uarchalfr}[2]{\uarcxhalfr{#1}{#2}{.4}}

\nc{\bdarcx}[3]{\draw[blue](#1,0)arc(180:90:#3) (#1+#3,#3)--(#2-#3,#3) (#2-#3,#3) arc(90:0:#3);}
\nc{\bdarc}[2]{\darcx{#1}{#2}{.4}}
\nc{\rduarcx}[3]{\draw[red](#1,0)arc(180:270:#3) (#1+#3,0-#3)--(#2-#3,0-#3) (#2-#3,0-#3) arc(270:360:#3);}
\nc{\rduarc}[2]{\uarcx{#1}{#2}{.4}}
\nc{\duarcx}[3]{\draw(#1,0)arc(180:270:#3) (#1+#3,0-#3)--(#2-#3,0-#3) (#2-#3,0-#3) arc(270:360:#3);}
\nc{\duarc}[2]{\uarcx{#1}{#2}{.4}}

\nc{\uv}[1]{\fill (#1,2)circle(.1);}
\nc{\lv}[1]{\fill (#1,0)circle(.1);}
\nc{\stline}[2]{\draw(#1,2)--(#2,0);}
\nc{\tlab}[2]{\draw(#1,2)node[above]{\tiny $#2$};}
\nc{\tudots}[1]{\draw(#1,2)node{$\cdots$};}
\nc{\tldots}[1]{\draw(#1,0)node{$\cdots$};}

\nc{\uvw}[1]{\fill[white] (#1,2)circle(.1);}
\nc{\huv}[1]{\fill (#1,1)circle(.1);}
\nc{\llv}[1]{\fill (#1,-2)circle(.1);}
\nc{\arcup}[2]{
\draw(#1,2)arc(180:270:.4) (#1+.4,1.6)--(#2-.4,1.6) (#2-.4,1.6) arc(270:360:.4);
}
\nc{\harcup}[2]{
\draw(#1,1)arc(180:270:.4) (#1+.4,.6)--(#2-.4,.6) (#2-.4,.6) arc(270:360:.4);
}
\nc{\arcdn}[2]{
\draw(#1,0)arc(180:90:.4) (#1+.4,.4)--(#2-.4,.4) (#2-.4,.4) arc(90:0:.4);
}
\nc{\cve}[2]{
\draw(#1,2) to [out=270,in=90] (#2,0);
}
\nc{\hcve}[2]{
\draw(#1,1) to [out=270,in=90] (#2,0);
}
\nc{\catarc}[3]{
\draw(#1,2)arc(180:270:#3) (#1+#3,2-#3)--(#2-#3,2-#3) (#2-#3,2-#3) arc(270:360:#3);
}

\nc{\arcr}[2]{
\draw[red](#1,0)arc(180:90:.4) (#1+.4,.4)--(#2-.4,.4) (#2-.4,.4) arc(90:0:.4);
}
\nc{\arcb}[2]{
\draw[blue](#1,2-2)arc(180:270:.4) (#1+.4,1.6-2)--(#2-.4,1.6-2) (#2-.4,1.6-2) arc(270:360:.4);
}
\nc{\loopr}[1]{
\draw[blue](#1,-2) edge [out=130,in=50,loop] ();
}
\nc{\loopb}[1]{
\draw[red](#1,-2) edge [out=180+130,in=180+50,loop] ();
}
\nc{\redto}[2]{\draw[red](#1,0)--(#2,0);}
\nc{\bluto}[2]{\draw[blue](#1,0)--(#2,0);}
\nc{\dotto}[2]{\draw[dotted](#1,0)--(#2,0);}
\nc{\lloopr}[2]{\draw[red](#1,0)arc(0:360:#2);}
\nc{\lloopb}[2]{\draw[blue](#1,0)arc(0:360:#2);}
\nc{\rloopr}[2]{\draw[red](#1,0)arc(-180:180:#2);}
\nc{\rloopb}[2]{\draw[blue](#1,0)arc(-180:180:#2);}
\nc{\uloopr}[2]{\draw[red](#1,0)arc(-270:270:#2);}
\nc{\uloopb}[2]{\draw[blue](#1,0)arc(-270:270:#2);}
\nc{\dloopr}[2]{\draw[red](#1,0)arc(-90:270:#2);}
\nc{\dloopb}[2]{\draw[blue](#1,0)arc(-90:270:#2);}
\nc{\llloopr}[2]{\draw[red](#1,0-2)arc(0:360:#2);}
\nc{\llloopb}[2]{\draw[blue](#1,0-2)arc(0:360:#2);}
\nc{\lrloopr}[2]{\draw[red](#1,0-2)arc(-180:180:#2);}
\nc{\lrloopb}[2]{\draw[blue](#1,0-2)arc(-180:180:#2);}
\nc{\ldloopr}[2]{\draw[red](#1,0-2)arc(-270:270:#2);}
\nc{\ldloopb}[2]{\draw[blue](#1,0-2)arc(-270:270:#2);}
\nc{\luloopr}[2]{\draw[red](#1,0-2)arc(-90:270:#2);}
\nc{\luloopb}[2]{\draw[blue](#1,0-2)arc(-90:270:#2);}

\nc{\larcb}[2]{
\draw[blue](#1,0-2)arc(180:90:.4) (#1+.4,.4-2)--(#2-.4,.4-2) (#2-.4,.4-2) arc(90:0:.4);
}
\nc{\larcr}[2]{
\draw[red](#1,2-2-2)arc(180:270:.4) (#1+.4,1.6-2-2)--(#2-.4,1.6-2-2) (#2-.4,1.6-2-2) arc(270:360:.4);
}

\rnc{\H}{\mathscr H}
\rnc{\L}{\mathscr L}
\nc{\R}{\mathscr R}
\nc{\D}{\mathcal D}
\nc{\J}{\mathscr D}

\nc{\ssim}{\mathrel{\raise0.25 ex\hbox{\oalign{$\approx$\crcr\noalign{\kern-0.84 ex}$\approx$}}}}
\nc{\POI}{\mathcal{POI}}
\nc{\wb}{\overline{w}}
\nc{\ub}{\overline{u}}
\nc{\vb}{\overline{v}}
\nc{\fb}{\overline{f}}
\nc{\gb}{\overline{g}}
\nc{\hb}{\overline{h}}
\nc{\pb}{\overline{p}}
\rnc{\sb}{\overline{s}}
\nc{\XR}{\pres{X}{R\,}}
\nc{\YQ}{\pres{Y}{Q}}
\nc{\ZP}{\pres{Z}{P\,}}
\nc{\XRone}{\pres{X_1}{R_1}}
\nc{\XRtwo}{\pres{X_2}{R_2}}
\nc{\XRthree}{\pres{X_1\cup X_2}{R_1\cup R_2\cup R_3}}
\nc{\er}{\eqref}
\nc{\larr}{\mathrel{\hspace{-0.35 ex}>\hspace{-1.1ex}-}\hspace{-0.35 ex}}
\nc{\rarr}{\mathrel{\hspace{-0.35 ex}-\hspace{-0.5ex}-\hspace{-2.3ex}>\hspace{-0.35 ex}}}
\nc{\lrarr}{\mathrel{\hspace{-0.35 ex}>\hspace{-1.1ex}-\hspace{-0.5ex}-\hspace{-2.3ex}>\hspace{-0.35 ex}}}
\nc{\nn}{\tag*{}}
\nc{\epfal}{\tag*{$\Box$}}
\nc{\tagd}[1]{\tag*{(#1)$'$}}
\nc{\ldb}{[\![}
\nc{\rdb}{]\!]}
\nc{\sm}{\setminus}
\nc{\I}{\mathcal I}
\nc{\InSn}{\I_n\setminus\S_n}
\nc{\dom}{\operatorname{dom}} \nc{\codom}{\operatorname{dom}}
\nc{\ojin}{1\leq j<i\leq n}
\nc{\eh}{\widehat{e}}
\nc{\wh}{\widehat{w}}
\nc{\uh}{\widehat{u}}
\nc{\vh}{\widehat{v}}
\nc{\sh}{\widehat{s}}
\nc{\fh}{\widehat{f}}
\nc{\textres}[1]{\text{\emph{#1}}}
\nc{\aand}{\emph{\ and \ }}
\nc{\iif}{\emph{\ if \ }}
\nc{\textlarr}{\ \larr\ }
\nc{\textrarr}{\ \rarr\ }
\nc{\textlrarr}{\ \lrarr\ }

\nc{\comma}{,\ }

\nc{\COMMA}{,\quad}
\nc{\TnSn}{\T_n\setminus\S_n} 
\nc{\TmSm}{\T_m\setminus\S_m} 
\nc{\TXSX}{\T_X\setminus\S_X} 
\rnc{\S}{\mathcal S}

\nc{\T}{\mathcal T} 
\nc{\A}{\mathscr A} 
\nc{\B}{\mathcal B} 
\rnc{\P}{\mathcal P} 
\nc{\K}{\mathcal K}
\nc{\PB}{\mathcal{PB}} 
\nc{\rank}{\operatorname{rank}}

\nc{\mtt}{\!\!\!\mt\!\!\!}

\nc{\sub}{\subseteq}
\nc{\la}{\langle}
\nc{\ra}{\rangle}
\nc{\mt}{\mapsto}
\nc{\im}{\mathrm{im}}
\nc{\id}{\mathrm{id}}
\nc{\bn}{\mathbf{n}}
\nc{\ba}{\mathbf{a}}
\nc{\bl}{\mathbf{l}}
\nc{\bm}{\mathbf{m}}
\nc{\bk}{\mathbf{k}}
\nc{\br}{\mathbf{r}}
\nc{\al}{\alpha}
\nc{\be}{\beta}
\nc{\ga}{\gamma}
\nc{\Ga}{\Gamma}
\nc{\de}{\delta}
\nc{\ka}{\kappa}
\nc{\lam}{\lambda}
\nc{\Lam}{\Lambda}
\nc{\si}{\sigma}
\nc{\Si}{\Sigma}
\nc{\oijn}{1\leq i<j\leq n}
\nc{\oijm}{1\leq i<j\leq m}

\nc{\comm}{\rightleftharpoons}
\nc{\AND}{\qquad\text{and}\qquad}

\nc{\bit}{\vspace{-3 truemm}\begin{itemize}}
\nc{\bmc}{\vspace{-3 truemm}\begin{itemize}\begin{multicols}}
\nc{\emc}{\end{multicols}\end{itemize}\vspace{-3 truemm}}
\nc{\eit}{\end{itemize}\vspace{-3 truemm}}
\nc{\ben}{\vspace{-3 truemm}\begin{enumerate}}
\nc{\een}{\end{enumerate}\vspace{-3 truemm}}
\nc{\eitres}{\end{itemize}}

\nc{\set}[2]{\{ {#1} : {#2} \}} 
\nc{\bigset}[2]{\big\{ {#1}: {#2} \big\}} 
\nc{\Bigset}[2]{\Big\{ \,{#1}\, \,\Big|\, \,{#2}\, \Big\}}

\nc{\pres}[2]{\la {#1} \,|\, {#2} \ra}
\nc{\bigpres}[2]{\big\la {#1} \,\big|\, {#2} \big\ra}
\nc{\Bigpres}[2]{\Big\la \,{#1}\, \,\Big|\, \,{#2}\, \Big\ra}
\nc{\Biggpres}[2]{\Bigg\la {#1} \,\Bigg|\, {#2} \Bigg\ra}

\nc{\pf}{\noindent{\bf Proof.}  }
\nc{\epf}{\hfill$\Box$\bigskip}
\nc{\epfres}{\hfill$\Box$}
\nc{\pfnb}{\pf}
\nc{\epfnb}{\bigskip}
\nc{\pfthm}[1]{\bigskip \noindent{\bf Proof of Theorem \ref{#1}}\,\,  } 
\nc{\pfprop}[1]{\bigskip \noindent{\bf Proof of Proposition \ref{#1}}\,\,  } 
\nc{\epfreseq}{\tag*{$\Box$}}

\makeatletter
\newcommand\footnoteref[1]{\protected@xdef\@thefnmark{\ref{#1}}\@footnotemark}
\makeatother

\numberwithin{equation}{section}

\newtheorem{thm}[equation]{Theorem}
\newtheorem{lemma}[equation]{Lemma}
\newtheorem{cor}[equation]{Corollary}
\newtheorem{prop}[equation]{Proposition}

\theoremstyle{definition}

\newtheorem{rem}[equation]{Remark}
\newtheorem{defn}[equation]{Definition}
\newtheorem{eg}[equation]{Example}
\newtheorem{ass}[equation]{Assumption}

\title{Variants of finite full transformation semigroups} 
\author{
Igor Dolinka%
\footnote{The first named author gratefully acknowledges the support of Grant No.~174019 of the Ministry of Education, Science, and Technological Development of the Republic of Serbia, and Grant No.~1136/2014 of the Secretariat of Science and Technological Development of the Autonomous Province of Vojvodina.}
\\
{\footnotesize \emph{Department of Mathematics and Informatics}}\\
{\footnotesize \emph{University of Novi Sad, Trg Dositeja Obradovi\'ca 4, 21101 Novi Sad, Serbia}}\\
{\footnotesize {\tt dockie@dmi.uns.ac.rs}}\\~\\
James East\\
{\footnotesize \emph{Centre for Research in Mathematics; School of Computing, Engineering and Mathematics}}\\
{\footnotesize \emph{University of Western Sydney, Locked Bag 1797, Penrith NSW 2751, Australia}}\\
{\footnotesize {\tt J.East\,@\,uws.edu.au}}
}

\date{}

\maketitle

\begin{abstract}
The \emph{variant} of a semigroup $S$ with respect to an element $a\in S$, denoted $S^a$, is the semigroup with underlying set $S$ and operation $\star$ defined by $x\star y=xay$ for $x,y\in S$.  In this article, we study variants $\TXa$ of the full transformation semigroup $\T_X$ on a finite set~$X$.  We explore the structure of $\TXa$ as well as its subsemigroups $\RegTXa$ (consisting of all regular elements) and $\EXa$ (consisting of all products of idempotents), and the ideals of $\RegTXa$.  Among other results, we calculate the rank and idempotent rank (if applicable) of each semigroup, and (where possible) the number of (idempotent) generating sets of the minimal possible size.


\end{abstract}

\section{Introduction}\label{sect:intro}

In John Howie's famous 1966 paper \cite{Howie1966}, it was shown that the semigroup $\Sing_X$ of all singular transformations on a finite set $X$ (i.e., all non-invertible functions $X\to X$) is generated by it idempotents.  In subsequent works, and with other authors, Howie calculated the rank (minimal size of a generating set) and idempotent rank (minimal size of an idempotent generating set) of $\Sing_X$ \cite{Gomes1987,Howie1978}; classified the idempotent generating sets of $\Sing_X$ of minimal size \cite{Howie1978}; calculated the rank and idempotent rank of the ideals of $\Sing_X$ \cite{Howie1990}; investigated the length function on $\Sing_X$ with respect to the generating set consisting of all idempotents of defect $1$ \cite{HLF}; and extended these results to various other kinds of transformation semigroups and generating sets \cite{GH1992,AAH2005,AAUH2008,Gomes1987,Howie7071}.  These works have been enormously influential, and have led to the development of several vibrant areas of research covering semigroups of (partial) transformations, matrices, partitions, endomorphisms of various algebraic structures, and more; see for example \cite{Fountain1992,Erdos1967,Fountain1991,EastGray,JEpnsn,Maltcev2007,Gray2007, Gray2008,Fountain1993,EF,AABK,SS2013} and references therein.  The current article continues in the spirit of this program of research, but takes it in a different direction; rather than concentrating on semigroups whose elements are variations of transformations of a set, we investigate semigroups of transformations under natural alternative binary operations, studying the so-called \emph{variants} of the full transformation semigroup.

The study of semigroup variants goes back to the 1960 monograph of Lyapin \cite{Lyapin} and a 1967 paper of Magill \cite{Magill1967} that considers semigroups of functions $X\to Y$ under an operation defined by $f\cdot g=f\circ \theta\circ g$, where $\theta$ is some fixed function $Y\to X$; see also \cite{Kemprasit2002,Symons1975b,Chase1979,Thornton1982,MS1975,MS1978,DE2015_1,Brown1955}.  In the case that $X=Y$, this provides an alternative product on the \emph{full transformation semigroup} $\T_X$ (consisting of all functions $X\to X$) that we will have more to say about below.  More generally, the \emph{variant} of a semigroup $S$ with respect to an element~$a\in S$ is the semigroup, denoted $S^a$, with underlying set $S$ and operation $\star$ defined by $x\star y=xay$ for each $x,y\in S$.  
Variants of arbitrary semigroups were first studied in 1983 by Hickey \cite{Hickey1983}, where (among other things) they were used to provide a novel characterisation of Nambooripad's celebrated partial order \cite{Nambooripad1980} on a regular semigroup; see also \cite{Hickey1986}.  As noted by Khan and Lawson \cite{KL2001}, variants arise naturally in relation to Rees matrix semigroups, and also provide a useful alternative to the group of units in some classes of non-monoidal regular semigroups (we explore the latter idea in Section \ref{sect:variants} below).

If $S$ is a group, it is easy to see that $S^a$ is isomorphic to $S$, the identity element of $S^a$ being $a^{-1}$; in a sense, this shows that no element of a group is more special than another, as the product may be ``scaled'' so that any element may play the role of the identity.  When $S$ is not a group, the situation can be very different.  Indeed, many semigroups with a relatively simple structure give rise to exceedingly complex variants; compare for example the right-most semigroup pictured in Figure \ref{fig:T1...T4} with some of its variants pictured in Figures~\ref{fig:V4_1233} and~\ref{fig:V4_1122_1222} (these figures are explained in detail below).\footnote{The authors are grateful to Attila Egri-Nagy for producing the GAP code for computing with semigroup variants.}  In complete contrast to the situation with groups, where every variant is isomorphic to the group itself, there exist semigroups for which all the variants are pairwise non-isomorphic; the bicyclic monoid is such a semigroup \cite{Tsyaputa2005}, and some more examples may be found in \cite{GMbook}.  On the other hand, some semigroups are isomorphic to all their variants (rectuangular bands, for example, which satisfy the identity $xay=xy$).

Variants of finite full transformation semigroups have been studied in a variety of contexts.  For example, Tsyaputa classified the non-isomorphic variants \cite{Tsyaputa2004}, characterised Green's relations \cite{Tsyaputa2005} and, together with Mazorchuk, classified the isolated subsemigroups \cite{MT2008}; see also \cite{Tsyaputa2005b,KT2005,Tsyaputa2006} where similar problems were considered in the context of partial transformations and partial permutations, and also \cite{Chase1979,Thornton1982,MS1975} where more general semigroups of functions and relations are considered.  The recent monograph of Ganyushkin and Mazorchuk \cite{GMbook} contains an entire chapter devoted to variants of various kinds of transformation semigroups, covering mostly Green's relations and the classification and enumeration of distinct variants.  In the current article, we take these existing results as our point of departure, and we investigate the kind of problems discussed in the opening paragraph in the context of the variants $\TXa$ of a finite full transformation semigroup~$\T_X$.
%
The structure and main results of the article are as follows.  In Sections~\ref{sect:TX} and~\ref{sect:variants}, we recall various facts regarding transformation semigroups and general variants (respectively), and also give a new characterisation of Green's relations on arbitrary variants (Proposition \ref{prop_green}); from these, we deduce Tsyaputa's above-mentioned results as corollaries in Section \ref{sect:TXa} (Theorem \ref{green_thm}), 
where we also explore the Green's structure of $\TXa$ further by investigating the natural partial order on the $\gD$-classes, using results regarding maximal $\gD$-classes to calculate the rank of $\TXa$ (Theorem \ref{thm_rankTXa}).  
%
%
The most substantial part of the article constitutes an investigation, in Section~\ref{sect:RegTXa}, of the structure of $\RegTXa$, the subsemigroup of $\TXa$ consisting of all regular elements (the elements of $\RegTXa$ are characterised in Section \ref{sect:TXa}, Proposition \ref{reg_prop}).  In particular, we identify $\RegTXa$ as a pullback product of the regular subsemigroups of two well-known semigroups consisting of transformations with restricted range and kernel (Propositions \ref{mono_prop} and \ref{epi_prop}), and we also show that $\RegTXa$ is a kind of ``inflation'' of the full transformation semigroup~$\T_A$, where $A$ denotes the image of $a$ (Theorem~\ref{inflation_thm}); among other things, these structural results allow us to calculate the size and rank of $\RegTXa$ (Corollary \ref{cor_|P|} and Theorem \ref{thm_rankP}).  The idempotent generated subsemigroup $\EXa$ of~$\TXa$ is studied in Section~\ref{sect:ETXa}, where we characterise the elements of $\EXa$ (Theorem \ref{thm_EXa}), calculate the rank and idempotent rank of $\EXa$ (showing in particular that these are equal, Theorem \ref{thm_rankEXa}), and classify and enumerate the minimal idempotent generating sets (Theorem \ref{thm_nrIGS}).  Finally, in Section \ref{sect:ideals}, we investigate the proper ideals of $\RegTXa$, showing that they are idempotent generated and calculating their rank and idempotent rank (which are again equal, Theorem \ref{thm_RegTXa_ideals}).

\section{Transformation semigroups}\label{sect:TX}

In this section, we record some basic notation and facts concerning finite transformation semigroups that we will need in what follows.  

If $S$ is any semigroup and $U\sub S$, we denote by $E(U)=\set{x\in U}{x^2=x}$ the set of all idempotents from $U$.  If $U\sub S$, we write $\la U\ra$ for the subsemigroup of $S$ generated by $U$, which consists of all products $u_1\cdots u_k$ where $k\geq1$ and $u_1,\ldots,u_k\in U$.  We write $\rank(S)$ for the \emph{rank of $S$}, defined to be the least cardinality of a subset $U\sub S$ such that $S=\la U\ra$.  If $S$ is idempotent generated, we write $\idrank(S)$ for the \emph{idempotent rank of $S$}, defined to be the least cardinality of a subset $U\sub E(S)$ such that $S=\la U\ra$.  Generation will always be in the variety of semigroups.

Recall that Green's relations $\gR$, $\gL$, $\gJ$, $\gH$, $\gD$, on a semigroup $S$ are defined, for $x,y\in S$, by
\begin{gather*}
x\gR y \iff x\Sone=y\Sone \COMMA x\gL y \iff \Sone x=\Sone y \COMMA x\gJ y \iff \Sone x\Sone =\Sone y\Sone  , \\
{\gH}={\gR\cap\gL} \COMMA {\gD}={\gR\circ\gL}={\gL\circ\gR}.
\end{gather*}
Here, $\Sone $ denotes the monoid obtained from $S$ by adjoining an identity element $1$, if necessary.  (We use the notation $\Sone$ rather than the more standard $S^1$ for reasons that will become clear shortly.)  If $x\in S$, and if $\gK$ is one of $\gR$, $\gL$, $\gJ$, $\gH$, $\gD$, we denote by $K_x$ the $\gK$-class of $x$ in~$S$.  An $\gH$-class contains an idempotent if and only if it is a group, in which case it is a maximal subgroup of $S$.  The $\gJ$-classes of $S$ are partially ordered; we say that $J_x\leq J_y$ if $x\in \Sone y\Sone$.  If $S$ is finite, then ${\gJ}={\gD}$.  An element $x\in S$ is \emph{regular} if $x=xyx$ and $y=yxy$ for some $y\in S$ or, equivalently, if $D_x$ contains an idempotent, in which case $R_x$ and $L_x$ do, too.  We write $\Reg(S)$ for the set of all regular elements of $S$, and we say $S$ is regular if $S=\Reg(S)$.

Let $X$ be a finite set with $|X|=n$.  The \emph{full transformation semigroup on $X$} is the (regular) semigroup $\T_X$ of all transformations of $X$ (i.e., all functions $X\to X$), under the operation of composition.  We write $xf$ for the image of $x\in X$ under $f\in\T_X$, and we compose functions from left to right.  If $f\in\T_X$, we will write
\[
f= \bigtrans{F_1 & \cdots & F_m}{f_1 & \cdots & f_m}
\]
to indicate that $X=F_1\sqcup\cdots\sqcup F_m$ and $F_if=f_i$ for each $i$.  (The symbol ``$\sqcup$'' denotes disjoint union.)  Usually this notation will imply that $f_1,\ldots,f_m$ are distinct, but occasionally this will not be the case, and we will always specify this.  
As usual, we denote the \emph{image}, \emph{kernel} and \emph{rank} of~$f\in\T_X$ by
\[
\im(f)=\set{xf}{x\in X} \COMMA  \ker(f)=\set{(x,y)\in X\times X}{xf=yf} \COMMA  \rank(f)=|\im(f)|=|X/\ker(f)|.
\]
We will sometimes write $\ker(f)=(F_1|\cdots|F_m)$ to indicate that $\ker(f)$ has equivalence classes $F_1,\ldots,F_m$, and this notation will always imply that the $F_i$ are pairwise disjoint and non-empty.  The \emph{symmetric group on $X$} is the set $\S_X$ of all permutations of $X$ (i.e., all invertible functions $X\to X$) and is the group of units of~$\T_X$.  
In the case that $X=\{1,\ldots,n\}$, we will write $\T_X=\T_n$ and $\S_X=\S_n$.  In general, if $k$ is a non-negative integer, we will write $\bk=\{1,\ldots,k\}$.  (So $\bk=\emptyset$ if $k=0$.)
A transformation $f\in\T_n$ will often be written as $f=[1f,\ldots,nf]$.
Green's relations on $\T_X$ are easy to describe; see for example \cite{Hig,Howie}.

\ms
\begin{prop}\label{prop_greenTX}
If $f\in\T_X$, where $X$ is a finite set with $|X|=n$, then
\bit
\itemit{i} $R_f =\set{g\in\T_X}{\ker(f)=\ker(g)}$,
\itemit{ii} $L_f =\set{g\in\T_X}{\im(f)=\im(g)}$,
\itemit{iii} $H_f =\set{g\in\T_X}{\ker(f)=\ker(g)\text{\emph{ and }}\im(f)=\im(g)}$,
\itemit{iv} $D_f =\set{g\in\T_X}{\rank(f)=\rank(g)}$. 
\eit
The $\gD$-classes of $\T_X$ form a chain: $D_1<\cdots<D_n$, where $D_m = \set{f\in\T_X}{\rank(f)=m}$ for each $m\in\bn$.  A group $\gH$-class contained in $D_m$ is isomorphic to $\S_m$. \epfres
\end{prop}

\begin{figure}[ht]
\begin{center}
\includegraphics[width=7.5mm]{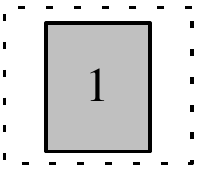}
\qquad
\includegraphics[width=10mm]{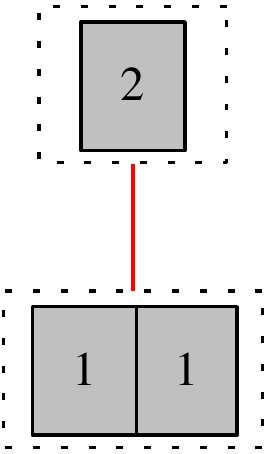}
\qquad
\includegraphics[width=14mm]{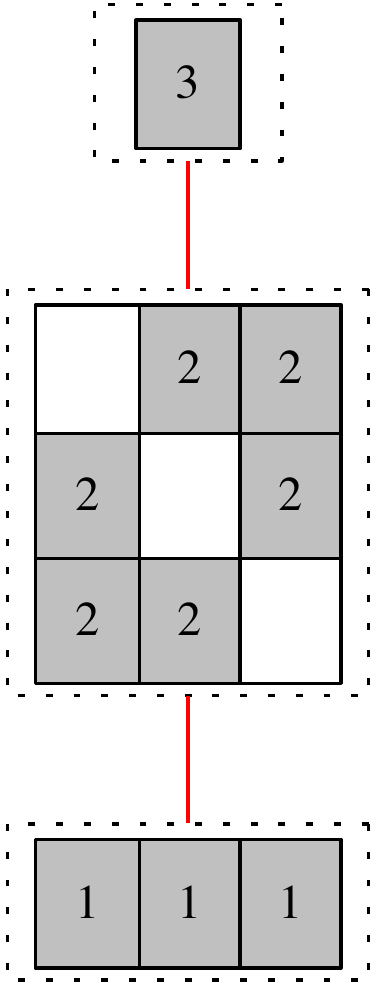}
\qquad
\includegraphics[width=25mm]{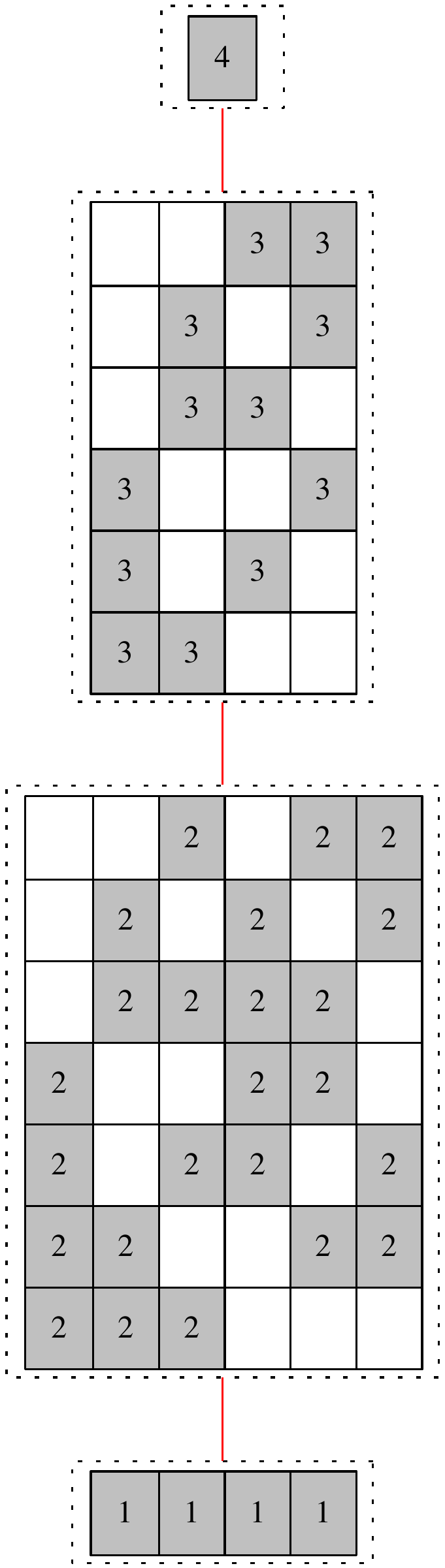}
    \caption{Egg box diagrams of the semigroups $\T_1$, $\T_2$, $\T_3$, $\T_4$ (left to right).}
    \label{fig:T1...T4}
   \end{center}
 \end{figure}

Note that $D_n=\S_X$.  For future reference, Figure \ref{fig:T1...T4} gives the so-called \emph{egg box diagrams} of the semigroups $\T_1$, $\T_2$, $\T_3$, $\T_4$.  Large boxes are $\gD$-classes; within a $\gD$-class, $\gR$-related (resp., $\gL$-related) elements are in the same row (resp., column); $\gH$-related elements are in the same cell; group $\gH$-classes are shaded grey and the label ``{\tt m}'' indicates that a given group is isomorphic to $\S_m$%
; the ${\gJ}={\gD}$-order is indicated by the edges between $\gD$-classes.  (See \cite{Hig,Howie} for more on egg box diagrams.)  The pictures were produced with the {\sc Semigroups} package \cite{GAP} on GAP \cite{GAP4}.

%
It is well-known that $\rank(\S_X)=2$ and $\rank(\T_X)=3$  if $|X|\geq3$; for example, $\S_n$ is generated by the transposition $[2,1,3,4,\ldots,n]$ and $n$-cycle $[2,3,4,\ldots,n,1]$, while $\T_X$ is generated by (any generating set for)~$\S_X$ along with any element of $D_{n-1}$; see for example \cite{Moo,Aiz,Vorobev1953,GMbook}.
The set $E(\T_X)$ of idempotents of $\T_X$ is not a subsemigroup, but the idempotent generated subsemigroup $\E_X=\la E(\T_X)\ra$ of $\T_X$ has a neat description.  For $x,y\in X$ with $x\not=y$, denote by $\ep_{xy}$ the (idempotent) transformation defined, for $z\in X$, by
\[
z\ep_{xy} = \begin{cases}
x &\text{if $z=y$}\\
z &\text{if $z\not=y$.}
\end{cases}
\]
Then $E(D_{n-1})=\set{\ep_{xy}}{x,y\in X,\ x\not=y}$.

\ms
\begin{thm}[Howie \cite{Howie1966,Howie1978}; Gomes and Howie \cite{Gomes1987}]\label{thm_TXSX}
If $X$ is a finite set with $|X|=n\geq2$, then
\[
\E_X=\la E(\T_X)\ra = \{1\}\cup(\TXSX) \AND \la E(D_{n-1}) \ra =\TXSX.
\]
Further, $\rank(\TXSX)=\idrank(\TXSX)=\rho_n$, where $\rho_2=2$ and $\rho_n={n\choose2}$ if $n\geq3$. \epfres
\end{thm}


%
%
%
The minimal idempotent generating sets of $\TXSX$ have a nice graphical interpretation.  Recall that a \emph{tournament on $X$} is a directed graph $\Ga$ with vertex set $X$ such that for each $x,y\in X$ with $x\not=y$, $\Ga$ contains precisely one of the directed edges $(x,y)$ or $(y,x)$.  Recall also that a directed graph on vertex set $X$ is \emph{strongly connected} if for any $x,y\in X$, there is a directed path from $x$ to $y$ in $\Ga$.  If $|X|\geq3$, we will write~$\bbT_X$ for the set of all strongly connected tournaments on $X$.  By convention, if $X=\{x,y\}$ is a set of size $2$, we will let $\bbT_X$ denote the set consisting of a single graph; namely, the graph with vertex set $X$ and directed edges $(x,y)$ and $(y,x)$.  For $U\sub E(D_{n-1})$, we define a graph $\Ga_U$ on vertex set $X$ with a directed edge $(x,y)$ corresponding to each $\ep_{xy}\in U$.


\ms
\begin{thm}[Howie \cite{Howie1978}]\label{thm_GaU}
Let $X$ be a finite set with $|X|=n\geq2$, and let
\[
U\sub E(D_{n-1}) = \set{\ep_{xy}}{x,y\in X,\ x\not=y}
\]
with $|U|=\rho_n$ (as defined in Theorem \ref{thm_TXSX}).  Then $\TXSX=\la U\ra$ if and only if $\Ga_U\in\bbT_X$.  In particular, the number of idempotent generating sets of the minimal size $\rho_n$ is equal to $|\bbT_X|$. \epfres
\end{thm}

\bs
\begin{rem}
A recurrence relation for the numbers $|\bbT_X|$ is given in \cite{Wright1970}.  The current authors have shown \cite{DE1} that any idempotent generating set for $\TXSX$ contains one of minimal possible size; a formula was also given for the total number of subsets of $E(D_{n-1})$ that generate $\TXSX$ (but are not necessarily of size $\rho_n$).  Arbitrary generating sets of minimal size were classified in \cite{AABK}.  The subsemigroup generated by the idempotents of an infinite transformation semigroup was described in \cite{Howie1966}.  
\end{rem}
%
%
%

\section{Variant semigroups}\label{sect:variants}

Let $S$ be a semigroup, and fix some element $a\in S$.  A new operation $\star_a$ may be defined on $S$ by
\[
x\star_a y =xay \qquad\text{for each $x,y\in S$.}
\]
We write $S^a$ for the semigroup $(S,\star_a)$ obtained in this fashion, and call $S^a$ the \emph{variant of $S$ with respect to~$a$}.  Since we fix $S$ and $a$ throughout this section, we will supress the subscript and simply write $\star$ for $\star_a$.  (Note that several authors write $\circ_a$ instead of $\star_a$, but we use the current notation so as not to interfere with the usual use of $\circ$ to denote composition of functions in $\T_X$.) 

If $u,v\in \Sone$, the map $x\mt vxu$ defines a homomorphism $S^{uav}\to S^a$.  If $S$ is a monoid with identity $1$, we write $G(S)$ for the group of units of $S$; that is,
\[
G(S) = \set{x\in S}{(\exists y\in S)\ xy=yx=1}.
\]
We have already noted that $G(\T_X)=\S_X$.
If $S$ is a monoid and $u,v\in G(S)$ are units, then the above map $S^{uav}\to S^a$ is invertible and, hence, an isomorphism.  As a special case, if $a$ is a unit, the maps $x\mt xa$ and $x\mt ax$ define isomorphisms $S^a\to S=S^1$.  
As a result, we will typically concern ourselves only with the case that $a$ is not a unit (although $S$ may in fact be a monoid), and call $S^a$ a \emph{non-trivial variant} in this case.  Our main objects of study are the (non-trivial) variants of a finite full transformation semigroup $\T_X$, but in this section we will prove some general results concerning arbitrary variants.

\begin{figure}[ht]
\begin{center}
\includegraphics[width=\textwidth]{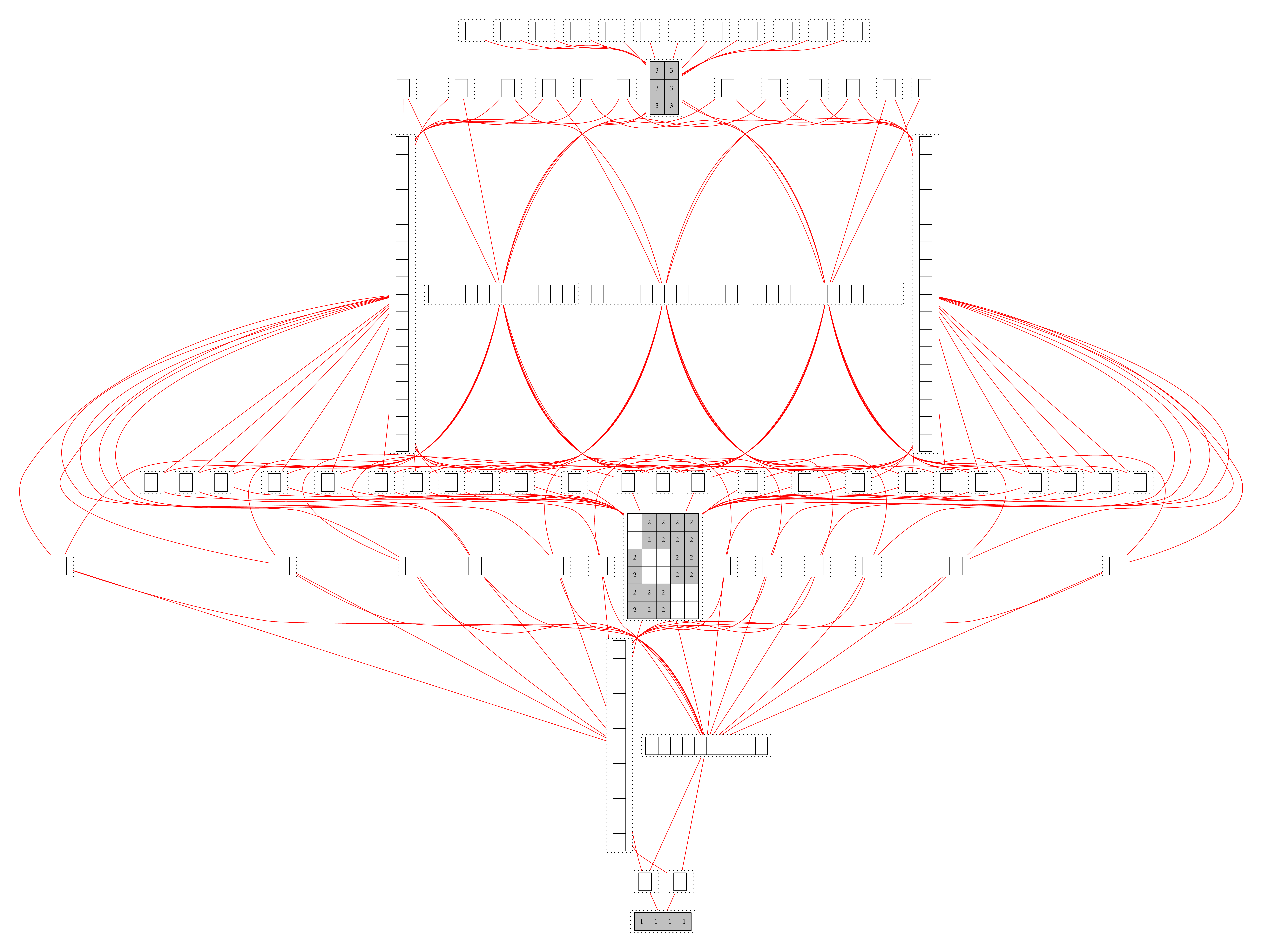}
    \caption{Egg box diagram of the  variant semigroup $\T_4^a$, where $a=[1,2,3,3]$.}
    \label{fig:V4_1233}
   \end{center}
 \end{figure}

\begin{figure}[ht]
\begin{center}
\includegraphics[width=\textwidth]{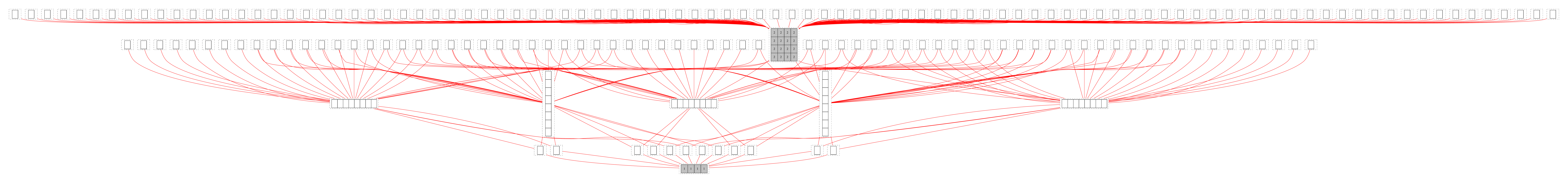}\\~\\
\includegraphics[width=\textwidth]{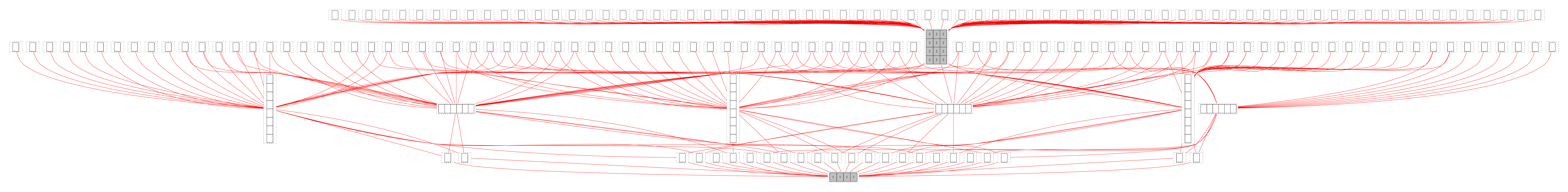}
    \caption{Egg box diagrams of the  variant semigroups $\T_4^a$, where $a=[1,1,2,2]$ (top) and $a=[1,2,2,2]$ (bottom).  The pdf file may be zoomed to obtain greater clarity.}
    \label{fig:V4_1122_1222}
   \end{center}
 \end{figure}

Before we do this, it is instructive to consider some examples; Figures \ref{fig:V4_1233} and \ref{fig:V4_1122_1222} illustrate the egg box diagrams of the variant semigroup $\T_4^a$ for various choices of $a\in\T_4$. 
%
%
A number of things become apparent when examining Figures \ref{fig:V4_1233} and \ref{fig:V4_1122_1222}.  In each case:
\bit
\item[(i)] $\T_4^a$ is not regular (as indicated by the many $\gD$-classes containing no idempotents).
\item[(ii)] A non-regular $\gD$-class of $\T_4^a$ is either a single $\gR$-class or a single $\gL$-class, or both (so a single $\gH$-class).
\item[(iii)] All the maximal $\gD$-classes are single $\gH$-classes (but a $\gD$-class consisting of a single $\gH$-class need not be maximal).
\item[(iv)] The number of maximal $\gD$-classes increases as $r=\rank(a)$ decreases.
\item[(v)] It is not evident from the picture, but every $\gH$-class contained in a non-regular $\gD$-class is a singleton.
\eit
In fact, all of these statements are true for arbitrary non-trivial variants of a finite transformation semigroup, while some are true for variants of arbitrary semigroups, as we will soon see.

We now prove a result concerning Green's relations on $S^a$.
In order to avoid confusion, if $\gK$ is one of $\gR$, $\gL$, $\gJ$, $\gH$, $\gD$, we will write $\gKa$ for Green's $\gK$-relation on the variant $S^a$, and write $K_x^a$ for the $\gKa$-class of $x\in S^a$.  It is easy to check that ${\gKa}\sub{\gK}$ for each relation $\gK$ and, hence, $K_x^a\sub K_x$ for each $x\in S$.
Throughout our investigations, a crucial role will be played by the sets
\[
P_1 = \set{x\in S}{xa \gR x} \COMMA P_2 = \set{x\in S}{ax \gL x}  \COMMA P=P_1\cap P_2.
\]
We note that $P_1=P_2=P=S$ if $S$ is a monoid and $a\in G(S)$ is a unit.

\ms 
\begin{lemma}\label{LyP1_RyP2}
If $y\in S$, then
\ss
\bmc2
\itemit{i} $y\in P_1$ if and only if $L_y\sub P_1$, 
\itemit{ii} $y\in P_2$ if and only if $R_y\sub P_2$.
\emc
The set $\Reg(S^a)$ of all regular elements of $S^a$ is contained in $P= P_1\cap P_2$.
\end{lemma}

\pf We just prove (i) because (ii) is dual.  
Suppose $y\in P_1$, and let $z\in L_y$ be arbitrary.  So $y\gR ya$, and we have $z=uy$ for some $u\in \Sone$.  But then $z=uy\gR uya=za$ since $\gR$ is a left congruence, so $z\in P_1$, whence $L_y\sub P_1$.  The other implication is trivial.  For the statement about regular elements, note that if $x\in\Reg(S^a)$, then $x=x\star y\star x=xayax$ for some $y\in S$. This gives $xa\gR x\gL ax$, so $x\in P$. \epf

\begin{prop}\label{prop_green}
If $x\in S$, then 
\ms
\bmc2
\itemit{i} $R_x^a = \begin{cases}
R_x\cap P_1 &\text{if $x\in P_1$}\\
\{x\} &\text{if $x\in S\sm P_1$,}
\end{cases}$
\itemit{ii} $L_x^a = \begin{cases}
L_x\cap P_2 &\hspace{0.8mm}\text{if $x\in P_2$}\\
\{x\} &\hspace{0.8mm}\text{if $x\in S\sm P_2$,}
\end{cases}
\phantom{
\begin{cases}
a\\b\\c\\d
\end{cases}
}
$
\itemit{iii} $H_x^a = \begin{cases}
H_x &\hspace{6.8mm}\text{if $x\in P$}\\
\{x\} &\hspace{6.8mm}\text{if $x\in S\sm P$,}
\end{cases}$
\itemit{iv} $D_x^a = \begin{cases}
D_x\cap P &\text{if $x\in P$}\\
L_x^a &\text{if $x\in P_2\sm P_1$}\\
R_x^a &\text{if $x\in P_1\sm P_2$}\\
\{x\} &\text{if $x\in S\sm (P_1\cup P_2)$.}
\end{cases}$
\emc
Further, if $x\in S\sm P$, then $H_x^a=\{x\}$ is a non-group $\gHa$-class of $S^a$.  
\end{prop}

\pf We begin with (i).  Suppose $y\in R_x^a\sm\{x\}$.  Then $x=y\star u=yau$ and $y=x\star v=xav$ for some $u,v\in S$.  But then $x=yau=xa(vau)$, so that $x\gR xa$, and $x\in P_1$.  In particular, if $x\in S\sm P_1$, then $R_x^a=\{x\}$.  Next, suppose $x\in P_1$.  If $y$ is another element of $R_x^a$ then, since $R_y^a=R_x^a$, the previous calculation shows that $y\in P_1$, and it follows that $R_x^a\sub P_1$.  Since we have already observed that $R_x^a\sub R_x$, it follows that $R_x^a\sub R_x\cap P_1$.  Conversely, suppose $y\in R_x\cap P_1$.  If $y=x$, then $y\in R_x^a$, so suppose $y\not=x$.  So $x=yu$ and $y=xv$ for some $u,v\in S$.  Also, $x=xaw$ and $y=yaz$ for some $w,z\in\Sone$, since $x,y\in P_1$.  Then $x=yu=yazu=y\star(zu)$ and, similarly, $y=x\star(wv)$, showing that $y\in R_x^a$.  

Part (ii) is dual to (i).  We now prove (iii).  If $x\in S\sm P$, then either $R_x^a=\{x\}$ or $L_x^a=\{x\}$ (or both).  In any case, $H_x^a=R_x^a\cap L_x^a=\{x\}$.  Next, suppose $x\in P$.  We have already noted that $H_x^a\sub H_x$.
Conversely, suppose $y\in H_x$.  If $y=x$, then $y\in H_x^a$, so suppose $y\not=x$.  Then $x=ys=ty$ and $y=xu=vx$ for some $s,t,u,v\in S$.  Also, $x=xaw=zax$ for some $w,z\in \Sone$, since $x\in P$.  But then $y=xu=xawu=x\star(wu)$ and $x=ys=vxs=vxaws=yaws=y\star(ws)$, showing that $y\gRa x$.  A similar calculation shows that $y\gLa x$, and we conclude that $y\in H_x^a$.

For part (iv), note that
\[
D_x^a = \bigcup_{y\in R_x^a}L_y^a  = \bigcup_{y\in L_x^a}R_y^a.
\]
In particular, if $x\in S\sm P_1$, then $R_x^a=\{x\}$, so that $D_x^a=L_x^a$.  Similarly, if $x\in S\sm P_2$, then $D_x^a=R_x^a$.  If $x\in S\sm(P_1\cup P_2)=(S\sm P_1)\cap(S\sm P_2)$, then $D_x^a=L_x^a=\{x\}$.  
Finally, if $x\in P$, then
\[
D_x^a = \bigcup_{y\in R_x^a}L_y^a = \bigcup_{y\in R_x\cap P_1}(L_y\cap P_2) = P_2\cap \bigcup_{y\in R_x\cap P_1} L_y = P_2\cap   \bigcup_{y\in R_x}(L_y\cap P_1 )  = P \cap \bigcup_{y\in R_x}L_y=P\cap D_x,
\]
where we have used parts (i) and (ii) in the second step, and Lemma \ref{LyP1_RyP2} (which tells us that $L_y\cap P_1$ is equal to $L_y$ or $\emptyset$ if $y\in P_1$ or $y\not\in P_1$, respectively) in the fourth step.

For the final statement about group $\gHa$-classes, suppose $H_x^a$ is a group, and let $e$ be the identity element of this group.  Then $x=x\star e=xae$ and also $x=eax$, so it follows that $xa\gR x\gL ax$, whence $x\in P$. \epf

\begin{rem}
In a sequel to the current paper \cite{DE2015_1}, we characterise the $\gJa$ 
relation, but we do not need this here. 
%
As noted above, if $S$ is a monoid and $a\in G(S)$ a unit, then $P_1=P_2=P=S$, in which case Green's relations on $S^a$ coincide exactly with the corresponding relations on $S\cong S^a$.  Let $x\in P=P_1\cap P_2$, and put $H=H_x^a=H_x$.  Whether $H$ is a group or non-group $\gH$-class of $S$ is independent of whether $H$ is a group or non-group $\gHa$-class of $S^a$.  See Table~\ref{tab:H_groups} for some examples with $S=\T_4$, $a=[1,2,3,3]$ and $x\in P$.  (See the next section for a description of the set $P$ in the case of $S=\T_X$.)
\end{rem}

\begin{table}[h]
\begin{center}
{
\begin{tabular}{|c|c|c|}
\hline
$x$ & Is $H_x$ a group $\gH$-class of $\T_4$? & Is $H_x$ a group $\gHa$-class of $\T_4^a$? \\
\hline
$[1,1,3,3]$ & Yes & Yes \\
$[4,2,2,4]$ & Yes & No \\
$[2,4,2,4]$ & No & Yes \\
$[1,3,1,3]$ & No & No \\
\hline
\end{tabular}
}
\end{center}
\caption{Group/non-group relationships between $H_x$ and $H_x^a$ in $\T_4$ and $\T_4^a$, where $a=[1,2,3,3]$.}
\label{tab:H_groups}
\end{table}

If $S$ is a monoid and $a\in G(S)$, then $S^a$ is a monoid (since then $S^a\cong S$).
The converse of this statement is also true, as we now demonstrate.  Part of the next proof is similar to that of \cite[Proposition 13.1.1]{GMbook}. 

\ms
\begin{prop}
Let $S$ be a semigroup and let $a\in S$.  Then $S^a$ is a monoid if and only if $S$ is a monoid and $a\in G(S)$, in which case $S^a$ is isomorphic to $S$.
\end{prop}

\pf It suffices to show the forwards implication, so suppose $S^a$ is a monoid with identity $e$.  In particular, for each $x\in S$, $x=x\star e=e\star x$; that is, $x=xae=eax$ for all $x$.  So $ae$ is a right identity for $S$, and $ea$ a left identity.  It follows that $ae=ea$ is a two sided identity for $S$, and that $a$ is a unit (with inverse $e$).  \epf

So $S^a$ is not a monoid in general, even if $S$ is itself a monoid.  The idea of the group of units of a monoid may be generalised to a non-unital semigroup $S$ by considering the so-called \emph{regularity presering} elements of $S$ \cite{KL2001,Hickey1983}; namely, those elements $a\in S$ for which $S^a$ is a regular semigroup.  The set of all regularity preserving elements of $S$ is denoted $\RP(S)$.  As the use of the word ``preserving'' suggests, $S$ can only contain regularity preserving elements if $S$ is itself regular, as may easily be checked (though there are regular semigroups $S$ for which $\RP(S)=\emptyset$, one example being $S=\TXSX$).  It is also clear that if $a\in\RP(S)$, then $J_a$ must be a maximum element in the ordering of $\gJ$-classes.  If $S$ is a regular monoid, then $\RP(S)=G(S)$, and this is just one of the reasons that $\RP(S)$ is considered to be a good analogue of the group of units in the case that $S$ is not a monoid.  The next result summarises some of the facts from \cite{KL2001} that we will need when investigating regularity preserving elements later.  
Recall that an element $u\in S$ is a \emph{mididentity} (sometimes called a \emph{midunit} or \emph{middle unit}) if $xuy=xy$ for all $x,y\in S$.  Semigroups with mididentity were first studied in \cite{Yamada1955} (the idea is also present in \cite{Thierrin1955}), and then more systematically in \cite{Ault1973,Ault1974}; the connection with semigroup variants is elucidated in \cite{KL2001,Hickey1983,BH1984}.

\ms
\begin{prop}[Khan and Lawson \cite{KL2001}]\label{RPS_prop}
Let $S$ be a regular semigroup.  
\bit
\itemit{i} An element $a\in S$ is regularity preserving if and only if $a\gH e$ for some regularity preserving idempotent $e\in E(S)$. (In particular, $\RP(S)$ is a union of groups.)
\itemit{ii} An idempotent $e\in E(S)$ is regularity preserving if and only if $fe \gR f \gL ef$ for all idempotents $f\in E(S)$.
\itemit{iii} Any mididentity is regularity preserving. \epfres
\eit
\end{prop}

So $S^a$ is not regular in general, even though $S$ may be regular itself.  But in some cases, $\Reg(S^a)$, the set of all regular elements of $S_a$, is a subsemigroup of $S$.  The next result was proved in \cite{KL2001} under the assumption that $S$ is regular, but the proof given there works unmodified for the following stronger statement.

\ms
\begin{lemma}[Khan and Lawson \cite{KL2001}]\label{lemma_RegSa}
Suppose $S$ is a semigroup, and that $aSa\sub\Reg(S)$ for some $a\in S$.  Then $\Reg(S^a)$ is a (regular) subsemigroup of $S^a$. \epfres
\end{lemma}

\section{The variant semigroup $\TXa$}\label{sect:TXa}

We now turn our attention to the main object of our study; namely, the variants $\TXa$, where $X$ is a finite set with $|X|=n$ and $a\in\T_X$.  The main results of this section include a characterisation of Green's relations and the ordering on ${\gJ}={\gD}$-classes, and the calculation of $\rank(\TXa)$.

It is easy to see that for any $a\in\T_X$, there is a permutation $p\in\S_X=G(\T_X)$ such that $ap\in E(\T_X)$ is an idempotent.  As noted in the previous section, $\TXa$ and $\T_X^{ap}$ are then isomorphic, so it suffices to assume that~$a$ is an idempotent.  
So for the remainder of the article, we fix an idempotent $a\in E(\T_X)$ with $r=\rank(a)$, and we write
\[
a = \left(\begin{matrix}
A_1 & \cdots & A_r \\
a_1 & \cdots & a_r
\end{matrix}\right).
\]
The condition that $a$ is an idempotent is equivalent to saying that $a_i\in A_i$ for each $i\in\br$.  Further, we will write $A=\im(a)=\{a_1,\ldots,a_r\}$ and $\al=\ker(a)=(A_1|\cdots|A_r)$.  We will also write $\lam_i=|A_i|$ for each~$i$, and for $I=\{i_1,\ldots,i_m\}\sub\br$, we define $\Lam_I=\lam_{i_1}\cdots\lam_{i_m}$.  In the special case that $I=\br$, we will write $\Lam=\Lam_{\br}=\lam_1\cdots\lam_r$.  As in the previous section, we will write $\star$ for $\star_a$.
If $r=n$, then $a\in\S_X=G(\T_X)$ and so, as we have noted, $\TXa\cong\T_X$.  All the problems we consider have been solved for $\T_X$, so we will assume throughout that $r<n$.  In particular, $\TXa$ is not a monoid, nor regular since $\S_X=\RP(\T_X)$.

As in the previous section, we will write $\gR$, $\gL$, $\gH$, ${\gD}={\gJ}$ for Green's relations on $\T_X$, and $\gRa$, $\gLa$, $\gHa$, ${\gDa}={\gJa}$ for Green's relations on $\TXa$.  If $f\in\T_X$ and if $\gK$ is one of $\gR,\gL,\gH,\gD$, we write $K_f$ and $K_f^a$ for the $\gK$-class and $\gKa$-class of $f$, respectively.  As we noted in the previous section for arbitrary variant semigroups, ${\gKa}\sub{\gK}$ for each $\gK$ and, hence, $K_f^a\sub K_f$ for each~$f$.  

As we have seen, the key to describing Green's relations on $\TXa$ are the sets
\[
P_1 = \set{f\in \T_X}{fa \gR f}\COMMA P_2 = \set{f\in \T_X}{af \gL f} \COMMA P=P_1\cap P_2.
\]
It will be convenient to have a more transparent characterisation of the elements of $P_1$ and $P_2$.  In order to give such a description, we introduce some terminology.  Let $B$ be a subset of $X$ and $\be$ an equivalence relation on $X$.  We say $B$ \emph{saturates} $\be$ if each $\be$-class contains at least one element of $B$.  We say $\be$ \emph{separates}~$B$ if each $\be$-class contains at most one element of $B$.  We call $B$ a \emph{cross-section} of $\be$ if $B$ saturates and is separated by $\be$.

\ms\ms
\begin{prop}\label{reg_prop}
\bit
\itemit{i} $P_1=\set{f\in\T_X}{\rank(fa)=\rank(f)}= \set{f\in \T_X}{\text{\emph{$\al$ separates $\im(f)$}}}$,
\itemit{ii} $P_2=\set{f\in\T_X}{\rank(af)=\rank(f)}=\set{f\in \T_X}{\text{\emph{$A$ saturates $\ker(f)$}}}$,
\itemit{iii} $P=\set{f\in\T_X}{\rank(afa)=\rank(f)}=\RegTXa$ is the set of all regular elements of $\TXa$, and is a subsemigroup of $\TXa$.
\eit
\end{prop}

\ms
\pf Let $f\in\T_X$ and write 
$
f = \trans{F_1 & \cdots & F_m}{f_1 & \cdots & f_m},
$
where $m=\rank(f)$.  For each $i\in\bm$, let $k_i\in\br$ be such that $f_i\in A_{k_i}$.  Note that
\[
f\in P_1 \iff fa\gR f \iff \ker(fa)=\ker(f) \iff \rank(fa)=\rank(f),
\]
since $X$ is finite.  Note that for each $i\in\bm$, $F_ifa=f_ia=a_{k_i}$.  It follows that $\rank(fa)=m$ if and only if the set $\{k_1,\ldots,k_m\}$ has cardinality $m$, and this is clearly equivalent to $\al$ separating $\im(f)$, establishing (i).  

A similar argument shows that $f\in P_2$ if and only if $\rank(af)=\rank(f)$.  Next, note that $\im(af)\sub\im(f)$ and that for all $i\in\bm$, $f_i(af)^{-1}=F_ia^{-1}=\bigcup_{a_j\in F_i}A_j$.  So $\rank(af)=m$ if and only if $F_i\cap A\not=\emptyset$ for all~$i$, and this is clearly equivalent to $A$ saturating $\ker(f)$, giving (ii).  

Combining the arguments of the previous two paragraphs shows that $f\in P=P_1\cap P_2$ if and only if $\rank(afa)=\rank(f)$.  We have already seen in Lemma \ref{LyP1_RyP2} that $\RegTXa\sub P$.  Conversely, suppose $f\in P$.  Since $\rank(fa)=\rank(af)=m$, we may write
$
fa = \trans{F_1 & \cdots & F_m}{a_{k_1} & \cdots & a_{k_m}}
$ and $
af = \trans{G_1 & \cdots & G_m}{f_1 & \cdots & f_m},
$
where ${k_1},\ldots,{k_m}$ are distinct, and $G_1,\ldots,G_m$ are non-empty and pairwise disjoint.  Let $g\in\T_X$ be any transformation for which $a_{k_i}g\in G_i$ for each $i\in\bm$.  Then clearly,
$f=(fa)g(af)=f\star g\star f$, showing that $f\in\RegTXa$.  Finally, Lemma~\ref{lemma_RegSa} tells us that $P$ is a subsemigroup of $\TXa$. \epf

Note that if $\rank(f)>r$, then $f$ belongs to neither $P_1$ nor $P_2$.  The next result follows from Proposition~\ref{prop_green}.  Together with Proposition \ref{reg_prop}, it yields the characterisation of Green's relations on $\TXa$ given by Tsyaputa~\cite{Tsyaputa2004}; see also \cite[Theorem 13.4.2]{GMbook}.  

\ms
\begin{thm}\label{green_thm}
If $f\in\TXa$, then   
\ms
\bmc2
\itemit{i} $R_f^a = \begin{cases}
R_f\cap P_1 &\text{if $f\in P_1$}\\
\{f\} &\text{if $f\in \T_X\sm P_1$,}
\end{cases}$
\itemit{ii} $L_f^a = \begin{cases}
L_f\cap P_2 &\hspace{0.7mm}\text{if $f\in P_2$}\\
\{f\} &\hspace{0.7mm}\text{if $f\in \T_X\sm P_2$,}
\end{cases}
\phantom{
\begin{cases}
a\\b\\c\\d
\end{cases}
}
$

\itemit{iii} $H_f^a = \begin{cases}
H_f &\hspace{6.8mm}\text{if $f\in P$}\\
\{f\} &\hspace{6.8mm}\text{if $f\in \T_X\sm P$,}
\end{cases}$
\itemit{iv} $D_f^a = \begin{cases}
D_f\cap P &\text{if $f\in P$}\\
L_f^a &\text{if $f\in P_2\sm P_1$}\\
R_f^a &\text{if $f\in P_1\sm P_2$}\\
\{f\} &\text{if $f\in \T_X\sm (P_1\cup P_2)$.}
\end{cases}$
\emc
The sets $P_1$ and $P_2$ are described in Proposition \ref{reg_prop}.  In particular, $R_f^a=L_f^a=H_f^a=D_f^a=\{f\}$ if $\rank(f)>r$.  If $f\in\T_X\sm P$, then $H_f^a=\{f\}$ is a non-group $\gHa$-class of $\TXa$.  \epfres 
\end{thm}

\ms
\begin{rem}
The article \cite{MS1975} characterises Green's relations and the regular elements of the more general semigroup $T(X,Y,a)$ consisting of all functions $f:X\to Y$ under the operation $f\cdot g=f\circ a\circ g$, where $a:Y\to X$ is some fixed function and $\circ$ denotes the usual composition of functions.  This characterisation is, by necessity, far more complex than that given in Proposition \ref{reg_prop} and Theorem \ref{green_thm}.
\end{rem}

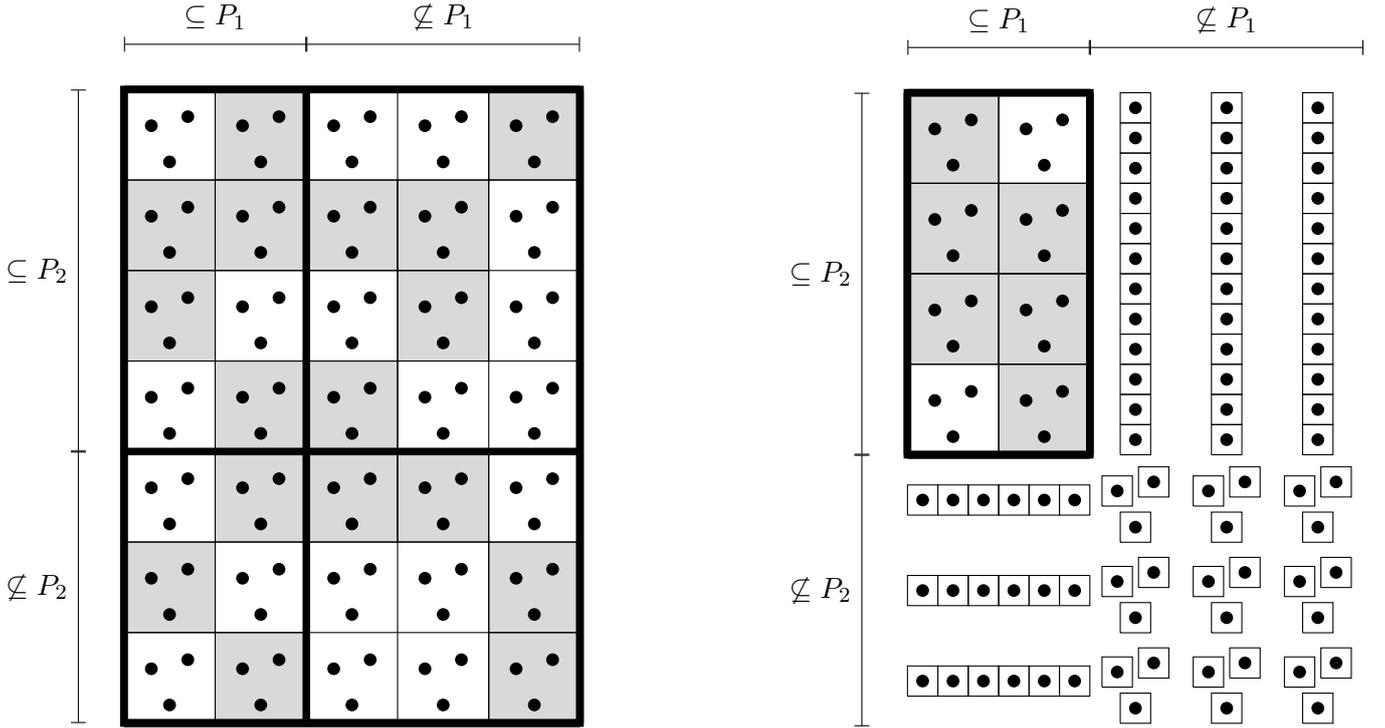
\begin{figure}[ht]
\begin{center}
\begin{tikzpicture}[scale=1.2]
\fillbox05
\fillbox04
\fillbox01
\fillbox16
\fillbox15
\fillbox13
\fillbox12
\fillbox10
\fillbox25
\fillbox23
\fillbox22
\fillbox35
\fillbox34
\fillbox32
\fillbox46
\fillbox41
\fillbox40
\foreach \x in {0,...,5}
{\draw (\x,0)--(\x,7);}
\foreach \x in {0,...,7}
{\draw (0,\x)--(5,\x);}
\foreach \x in {0,2,5}
{\draw[line width=1mm] (\x,0)--(\x,7);}
\foreach \x in {0,3,7}
{\draw[line width=1mm] (0,\x)--(5,\x);}
\draw[line width=1mm] (0,0)--(5,0)--(5,7)--(0,7)--(0,0)--(5,0);
\draw[|-|] (0,7.5)--(2,7.5);
\draw[|-|] (2,7.5)--(5,7.5);
\draw[|-|] (-.5,0)--(-.5,3);
\draw[|-|] (-.5,7)--(-.5,3);
\draw(1,7.8)node{$\sub P_1\phantom{}$};
\draw(3.5,7.8)node{$\not\sub P_1\phantom{}$};
\draw(-.5,5)node[left]{$\sub P_2\phantom{}$};
\draw(-.5,1.5)node[left]{$\not\sub P_2\phantom{}$};
\foreach \x in {0,...,4}
\foreach \y in {0,...,6}
{\fill (\x+0.3,\y+0.6)circle(.07); \fill (\x+0.7,\y+0.7)circle(.07); \fill (\x+0.5,\y+0.2)circle(.07);}
\end{tikzpicture}
\hfill
\begin{tikzpicture}[scale=1.2]
\fillbox05
\fillbox04
\fillbox06
\fillbox15
\fillbox14
\fillbox13
\foreach \x in {0,...,2}
{\draw (\x,3)--(\x,7);}
\foreach \x in {3,...,7}
{\draw (0,\x)--(2,\x);}
\foreach \x in {0,2}
{\draw[line width=1mm] (\x,3)--(\x,7);}
\foreach \x in {3,7}
{\draw[line width=1mm] (0,\x)--(2,\x);}
\draw[line width=1mm] (0,3)--(2,3)--(2,7)--(0,7)--(0,3)--(2,3);
\draw[|-|] (0,7.5)--(2,7.5);
\draw[|-|] (2,7.5)--(5,7.5);
\draw[|-|] (-.5,0)--(-.5,3);
\draw[|-|] (-.5,7)--(-.5,3);
\draw(1,7.8)node{$\sub P_1\phantom{}$};
\draw(3.5,7.8)node{$\not\sub P_1\phantom{}$};
\draw(-.5,5)node[left]{$\sub P_2\phantom{}$};
\draw(-.5,1.5)node[left]{$\not\sub P_2\phantom{}$};
\foreach \x in {0,1}
\foreach \y in {3,4,5,6}
{\fill (\x+0.3,\y+0.6)circle(.07); \fill (\x+0.7,\y+0.7)circle(.07); \fill (\x+0.5,\y+0.2)circle(.07);}
\foreach \x in {2,3,4}
\foreach \y in {0,1,2}
{\fill (\x+0.3,\y+0.6)circle(.07); \fill (\x+0.7,\y+0.7)circle(.07); \fill (\x+0.5,\y+0.2)circle(.07);}
\foreach \x in {2,3,4}
\foreach \y in {0,1,2}
\foreach \z in {0.16666}
{ 
\draw (\x+0.3-\z,\y+0.6-\z)--(\x+0.3+\z,\y+0.6-\z)--(\x+0.3+\z,\y+0.6+\z)--(\x+0.3-\z,\y+0.6+\z)--(\x+0.3-\z,\y+0.6-\z);
\draw (\x+0.7-\z,\y+0.7-\z)--(\x+0.7+\z,\y+0.7-\z)--(\x+0.7+\z,\y+0.7+\z)--(\x+0.7-\z,\y+0.7+\z)--(\x+0.7-\z,\y+0.7-\z);
\draw (\x+0.5-\z,\y+0.2-\z)--(\x+0.5+\z,\y+0.2-\z)--(\x+0.5+\z,\y+0.2+\z)--(\x+0.5-\z,\y+0.2+\z)--(\x+0.5-\z,\y+0.2-\z);
}
\foreach \x in {.16666,.5,.833333,1.16666,1.5,1.83333}
\foreach \y in {0.5,1.5,2.5}
\foreach \z in {0.166666}
{ 
\fill (\x,\y)circle(.07);
\draw (\x-\z,\y-\z)--(\x+\z,\y-\z)--(\x+\z,\y+\z)--(\x-\z,\y+\z)--(\x-\z,\y-\z);
}
\foreach \x in {2.5,3.5,4.5}
\foreach \y in {3.16666,3.5,3.833333,4.16666,4.5,4.833333,5.16666,5.5,5.833333,6.16666,6.5,6.833333}
\foreach \z in {0.16666}
{ 
\fill (\x,\y)circle(.07);
\draw (\x-\z,\y-\z)--(\x+\z,\y-\z)--(\x+\z,\y+\z)--(\x-\z,\y+\z)--(\x-\z,\y-\z);
}
\end{tikzpicture}
    \caption{A schematic diagram of the way a $\gD$-class $D_m$ of $\T_X$ (with $2\leq m\leq r$) breaks up into $\gDa$-classes in $\TXa$.  Group $\gH$- and $\gHa$-classes are shaded grey.}
    \label{fig:green}
   \end{center}
 \end{figure}

Theorem \ref{green_thm} yields an intuitive picture of the Green's structure of $\TXa$.  Recall that the $\gD$-classes of $\T_X$ are precisely the sets $D_m=\set{f\in\T_X}{\rank(f)=m}$ for $1\leq m\leq n=|X|$.  Each of the $\gD$-classes $D_{r+1},\ldots,D_n$ separates completely into singleton $\gDa$-classes in $\TXa$.  (We will study these classes in more detail shortly.)  Next, note that $D_1\sub P$ (as the constant maps clearly belong to both $P_1$ and $P_2$), so $D_1$ remains a (regular) $\gDa$-class of $\TXa$.
Now fix some $2\leq m\leq r$, and recall that we are assuming that $r<n$.  The $\gD$-class $D_m$ is split into a single regular $\gDa$-class, namely $D_m\cap P$, and a number of non-regular $\gDa$-classes.  Some of these non-regular $\gDa$-classes are singletons, namely those of the form $D_f^a=\{f\}$ where $f\in D_m$ belongs to neither $P_1$ nor $P_2$.  Some of the non-regular $\gDa$-classes consist of one non-singleton $\gLa$-class, namely those of the form $D_f^a=L_f^a=L_f\cap P_2$, where $f\in D_m$ belongs to $P_2\sm P_1$; the $\gHa$-classes contained in such a $\gDa$-class are all singletons.  The remaining non-regular $\gDa$-classes in $D_m$ consist of one non-singleton $\gRa$-class, namely those of the form $D_f^a=R_f^a=R_f\cap P_1$, where $f\in D_m$ belongs to $P_1\sm P_2$; the $\gHa$-classes contained in such a $\gDa$-class are all singletons.  This is all pictured (schematically) in Figure \ref{fig:green}; see also Figures \ref{fig:V4_1233} and \ref{fig:V4_1122_1222}.

We now give some information about the order on the ${\gJa}={\gDa}$-classes of $\TXa$.  
Recall that in $\T_X$, $D_f\leq D_g$ if and only if $\rank(f)\leq \rank(g)$.  The situation is more complicated in $\TXa$.

\ms 
\begin{prop}\label{prop_maximalD}
Let $f,g\in\T_X$.  Then $D_f^a\leq D_g^a$ in $\TXa$ if and only if one of the following holds:
\ms
\bmc2
\itemit{i} $f=g$,
\itemit{ii} $\rank(f)\leq\rank(aga)$,
\itemit{iii} $\im(f)\sub\im(ag)$,
\itemit{iv} $\ker(f)\supseteq\ker(ga)$.
\emc
The maximal $\gDa$-classes are those of the form $D_f^a=\{f\}$ where $\rank(f)>r$.
\end{prop}

\pf Note that $D_f^a\leq D_g^a$ if and only if one of the following holds:
\ms
\bmc2
\item[(a)] $f=g$,
\item[(b)] $f=uagav$ for some $u,v\in\T_X$,
\item[(c)] $f=uag$ for some $u\in \T_X$,
\item[(d)] $f=gav$ for some $v\in\T_X$.
\emc
We clearly have the implications (b) $\Rightarrow$ (ii), (c) $\Rightarrow$ (iii), and (d) $\Rightarrow$ (iv).  Next, note that (ii) implies $D_f\leq D_{aga}$ in $\T_X$, from which (b) follows.  Next suppose (iii) holds.  Since $\im(f)\sub\im(ag)$, we may write
$
f= \trans{F_1 & \cdots & F_m}{f_1 & \cdots & f_m}
$ and $
ag = \trans{G_1 & \cdots & G_m & G_{m+1} & \cdots & G_l}{f_1 & \cdots & f_m & g_{m+1} & \cdots & g_l}.
$
For $i\in\bm$, let $g_i\in G_i$.  We then have $f=uag$, where
$
u = \trans{F_1 & \cdots & F_m}{g_1 & \cdots & g_m},
$
giving (c).  Finally, suppose (iv) holds, and write
$
f= \trans{F_1 & \cdots & F_m}{f_1 & \cdots & f_m}
$ and $
ga = \trans{G_1 & \cdots & G_l}{g_1 & \cdots & g_l}.
$
Since $\ker(f)\supseteq\ker(ga)$, there is a surjective function $q:\bl\to\bm$ such that $G_i\sub F_{iq}$ for all $i$.  We see then that $f=gav$, where $v\in\T_X$ is any transformation that extends the partial map
$
\trans{g_1 & \cdots & g_l}{f_{1q} & \cdots & f_{lq}},
$
giving (d).  

To prove the statement concerning maximal $\gDa$-classes, let $f\in\T_X$.  If $\rank(f)\leq r$, then $\rank(f)\leq \rank(a)=\rank(a1a)$, so that $D_f^a< D_1^a=\{1\}$, whence $D_f^a$ is not maximal.  (Here, $1\in\T_X$ denotes the identity element of $\T_X$, namely the identity map $X\to X$.)  On the other hand, suppose $\rank(f)>r$ and that $D_f^a\leq D_g^a$.  Then (ii) does not hold, since $\rank(aga)\leq\rank(a)=r<\rank(f)$.  Similarly, $\rank(ag)<\rank(f)$ and $\rank(ga)<\rank(f)$, so neither (iii) nor (iv) holds.  Having eliminated (ii--iv), we deduce that (i) must hold; that is, $f=g$, so $D_f^a=\{f\}$ is maximal.  \epf

\begin{figure}[ht]
\begin{center}
\includegraphics[width=\textwidth]{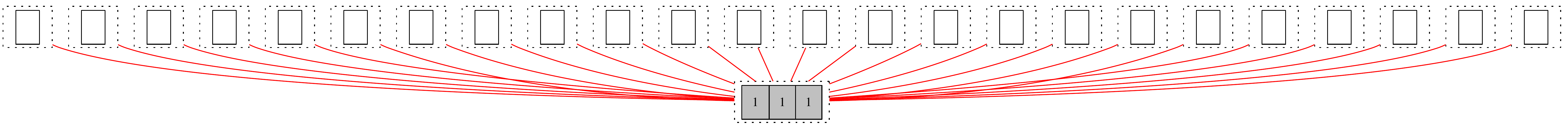}
    \caption{Egg box diagram of the  variant semigroup $\T_3^a$, where $a=[1,1,1]$.}
    \label{fig:V3_111}
   \end{center}
\end{figure}

\begin{rem}\label{rem_r=1}
If $r=\rank(a)=1$, then $\TXa$ has a very simple structure, as may be deduced from Theorem~\ref{green_thm} and Proposition \ref{prop_maximalD}; see Figure \ref{fig:V3_111} for an illustration in the case $n=|X|=3$.  This structure may also be observed directly.  For $x\in X$, denote by $c_x\in\T_X$ the constant map with image $\{x\}$.  If $a=c_x$, then for all $f,g\in\T_X$, $f\star g=fc_xg=c_xg=c_{xg}$.  
\end{rem}

The description of the maximal $\gDa$-classes from Proposition \ref{prop_maximalD} allows us to obtain information about $\rank(\TXa)$.  In order to avoid confusion when discussing generation, if $U\sub\T_X$, we will write $\la U\ra$ (resp.,~$\la U\raa$) for the subsemigroup of $\T_X$ (resp., $\TXa$) generated by $U$, which consists of all products $u_1\cdots u_k$ (resp., $u_1\star\cdots\star u_k$), where $k\geq1$ and $u_1,\ldots,u_k\in U$.  

\bs
\begin{thm}\label{thm_rankTXa}
Let $M=\set{f\in\T_X}{\rank(f)>r}$.  Then $\TXa=\la M\raa$.  Further, any generating set for $\TXa$ contains $M$.  Consequently, $M$ is the unique minimal (with respect to containment or size) generating set of~$\TXa$, and
\[
\rank(\TXa)=|M|=\sum_{m=r+1}^n S(n,m){n\choose m}m!,
\]
where $S(n,m)$ denotes the (unsigned) Stirling number of the second kind.
\end{thm}

\pf Consider the statement:
\nss
\begin{quote}
$H(m)$: ~ $\la M\raa$ contains $D_m\cup\cdots\cup D_n=\set{f\in\T_X}{\rank(f)\geq m}$.
\end{quote}
\nss
Since $H(1)$ says that $\TXa=\la M\raa$, it suffices to show that $H(m)$ is true for all $m\in\bn$.  We do this by (reverse) induction on $m$.  Note that $M=D_{r+1}\cup\cdots\cup D_n$, so $H(m)$ is clearly true for $m\geq r+1$.  Now suppose $H(m+1)$ is true for some $1\leq m\leq r$.  Let $f\in D_m$, and write
$
f= \trans{F_1 & \cdots & F_m}{f_1 & \cdots & f_m}.
$
Since $m\leq r<n$, we may assume that $|F_1|\geq2$.  Choose some non-trivial partition $F_1=F_1'\sqcup F_1''$.  Without loss of generality, we may also assume that $|A_1|\geq2$.  Choose some $a_1'\in A_1\sm\{a_1\}$, and put
$
g = \trans{F_1' & F_1'' & F_2 & \cdots & F_m}{a_1' & a_1 & a_2 & \cdots & a_m}.
$
So $g\in\la M\raa$, by the induction hypothesis.  Also, let $h\in\S_X\sub M$ be any permutation that extends the partial map
$
\trans{a_1 & \cdots & a_m}{f_1 & \cdots & f_m}.
$
Then $f=gah=g\star h\in\la M\raa$, so $H(m)$ is true, completing the inductive step.  

Any $f\in M$ belongs to a non-group, maximal $\gDa$-class, so it follows that any generating set of $\TXa$ must contain $M$.  This tells us that $M$ is the minimal generating set with respect to both size and containment, and that $\rank(\TXa)=|M|$.  The formula for $|M|$ follows from the well-known fact that
$
|D_m| = S(n,m){n\choose m}m!
$
for any $m\in \bn$ \cite{GMbook}.  This completes the proof. \epf

\begin{rem} 
It seems noteworthy that $\rank(\TXa)$ depends only on $r=\rank(a)$, and not on the sizes $\lam_1,\ldots,\lam_r$ of the kernel-classes of $a$.  See also Theorems \ref{thm_rankP}, \ref{thm_rankEXa} and \ref{thm_RegTXa_ideals}.
%
\end{rem}

The description of the order on $\gDa$-classes of $\TXa$ from Proposition \ref{prop_maximalD} may be simplified in the case that one of $f,g$ is regular.

\ms
\begin{prop}\label{prop_maximalD2}
Let $f,g\in\T_X$.
\bit
\itemit{i} If $f\in P$, then $D_f^a\leq D_g^a$ if and only if $\rank(f)\leq\rank(aga)$.
\itemit{ii} If $g\in P$, then $D_f^a\leq D_g^a$ if and only if $\rank(f)\leq\rank(g)$.
\eit
The regular $\gDa$-classes of $\TXa$ form a chain: $D_1^a<\cdots<D_r^a$, where $D_m^a=\set{f\in P}{\rank(f)=m}$ for $m\in\br$.
\end{prop}

\pf As in the proof of Proposition \ref{prop_maximalD}, $D_f^a\leq D_g^a$ if and only if one of the following holds:
\ms
\bmc2
\item[(a)] $f=g$,
\item[(b)] $f=uagav$ for some $u,v\in\T_X$,
\item[(c)] $f=uag$ for some $u\in \T_X$,
\item[(d)] $f=gav$ for some $v\in\T_X$.
\emc
Suppose first that $f\in P$, so $f=fahaf$ for some $h\in\T_X$.  Then (a) implies $f=fah(aga)haf$, (c) implies $f=u(aga)haf$, and (d) implies $f=fah(aga)v$.  So, in each of cases (a--d), we deduce that $\rank(f)\leq\rank(aga)$.  We have already observed that $\rank(f)\leq\rank(aga)$ implies $D_f^a\leq D_g^a$.  

Next, suppose $g\in P$.  Since $\rank(ag)=\rank(ga)=\rank(aga)=\rank(g)$, each of (a--d) implies $\rank(f)\leq\rank(g)$.  If $\rank(f)\leq\rank(g)=\rank(aga)$, then we already know that $D_f^a\leq D_g^a$.  The statement about regular $\gDa$-classes follows quickly from (ii).  \epf

Proposition \ref{prop_maximalD2} gives us some more partial information about the location of the ``fragmented'' $\gDa$-classes (see Figure \ref{fig:green}).  Specifically, a non-regular $\gDa$-class $D_f^a$ with $\rank(f)=m\leq r$ sits below $D_m^a$.  However, $D_f^a$ may or may not sit above $D_{m-1}^a$; this depends on $\rank(afa)$.  
For example, if $a=[1,1,1,4,5]$ and $f=[1,2,3,1,1]$, then $D_f^a$ sits between $D_1^a$ and $D_3^a$ but not above $D_2^a$ in $\T_5^a$.
%
While it is extremely difficult to enumerate all $\gDa$-classes (even maximal ones) that sit above $D_m^a$ but not $D_{m+1}^a$, where $m\in\br$ is arbitrary, we can enumerate those that sit right at the top of the picture, above the highest regular $\gDa$-class, $D_r^a$.  Recall that $\Lam=\lam_1\cdots\lam_r$, where $\lam_i=|A_i|$.

\ms
\begin{prop}
A maximal $\gDa$-class $D_f^a=\{f\}$ sits above $D_r^a$ in the ordering of $\gDa$-classes in~$\TXa$ if and only if $\rank(afa)=r<\rank(f)$.  The number of such $\gDa$-classes is equal to
$(n^{n-r}-r^{n-r})r!\Lam$.
\end{prop}

\pf The first statement follows from Proposition \ref{prop_maximalD2}(i).  It remains to count the number of transformations $f\in\T_X$ satisfying $\rank(afa)=r<\rank(f)$.  Note that such an $f$ maps $A$ to a cross-section of $\al=\ker(a)$.  The number of cross-sections of $\al$ is $\lam_1\cdots\lam_r=\Lam$, and once such a cross-section $B=\{b_1,\ldots,b_r\}$ is chosen, there are $r!$ ways to choose $f|_A$ (which maps $A$ bijectively to $B$).  There are $n^{n-r}-r^{n-r}$ ways to extend $f|_A$ to $f\in\T_X$ with $\rank(f)>r$. \epf

\section{The regular semigroup $\RegTXa$}\label{sect:RegTXa}

In this section, we study the subsemigroup
\[
P=\RegTXa=\set{f\in\T_X}{\rank(afa)=\rank(f)},
\]
consisting of all regular elements of~$\TXa$.  Key results include a description of $P$ as a subdirect product of the well-known semigroups $\RegTXA$ and $\RegTXal$ (see below for definitions), a realisation of $P$ as a kind of ``inflation'' of $\T_A$, combinatorial results on the number of Green's classes of certain types, and calculations of $|P|$ and $\rank(P)$.
As before, we assume that
\[
a = \left(\begin{matrix}
A_1 & \cdots & A_r \\
a_1 & \cdots & a_r
\end{matrix}\right)
\]
is an idempotent with $\rank(a)=r<n$, and we continue to write $A=\im(a)$, $\al=\ker(a)$, $\lam_i=|A_i|$, and so on.  
By Theorem~\ref{green_thm}, we see that $\RegTXa=D_1$ is a right zero semigroup in the case $r=1$ (see also Remark~\ref{rem_r=1} and Figure \ref{fig:V3_111}), in which case, all the problems we consider become trivial.  So for the duration of this section, we will assume that $1<r<n$.  

\begin{figure}[ht]
\begin{center}
\includegraphics[height=5.4mm]{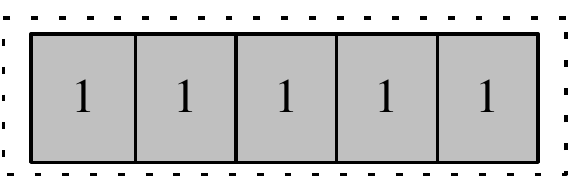}
\ 
\includegraphics[height=40mm]{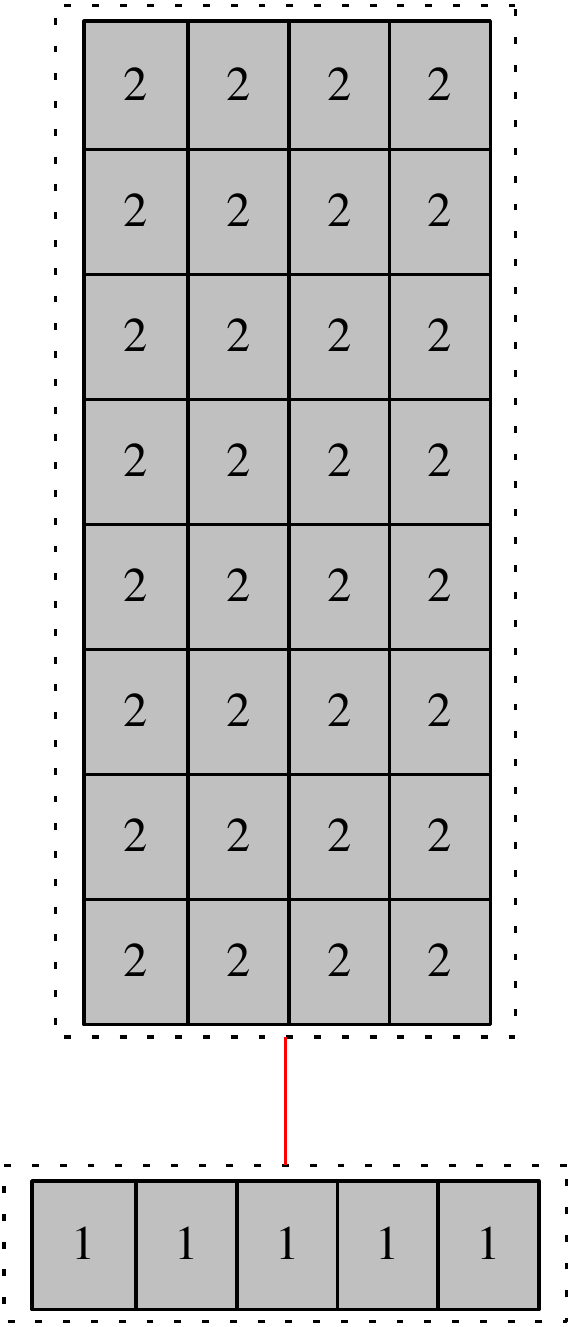}
\ 
\includegraphics[height=40mm]{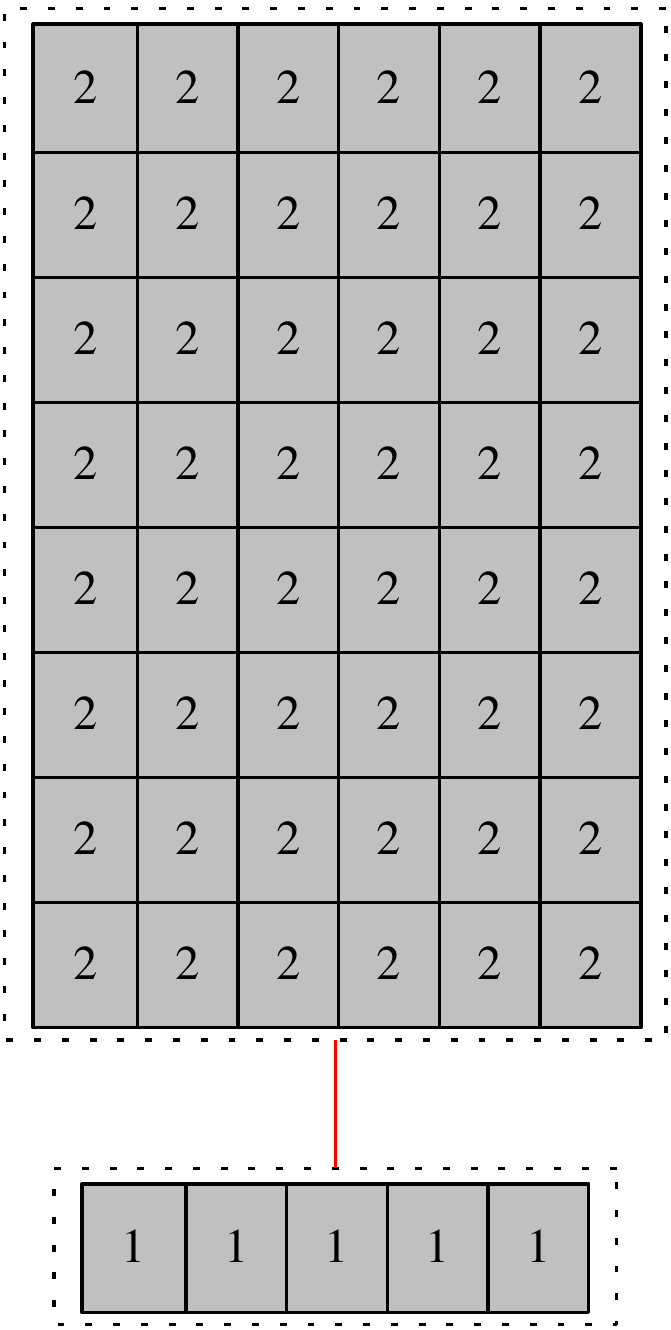}
\ 
\includegraphics[height=92mm]{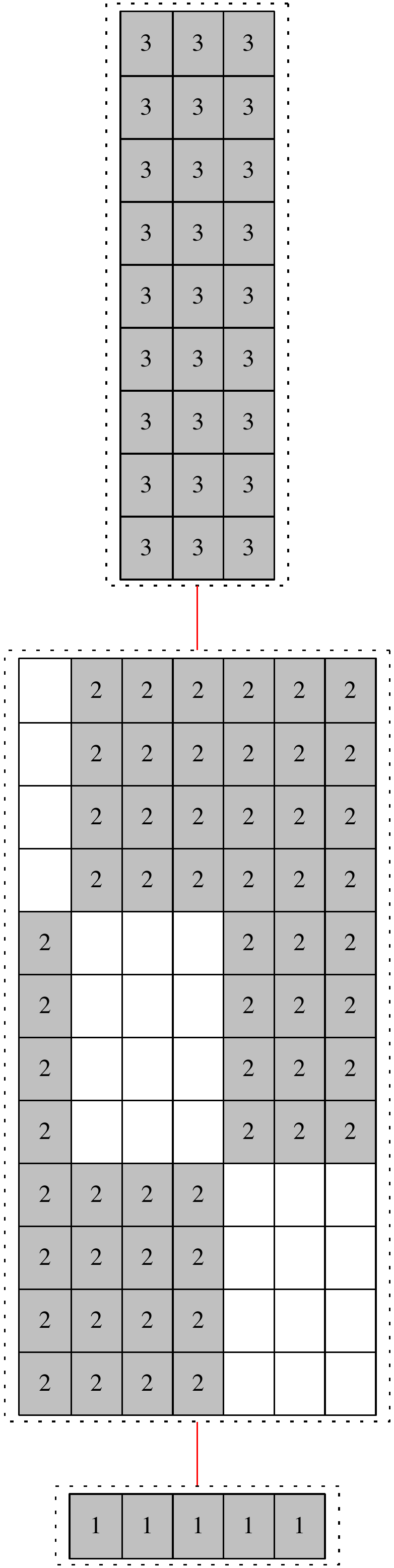}
\ 
\includegraphics[height=92mm]{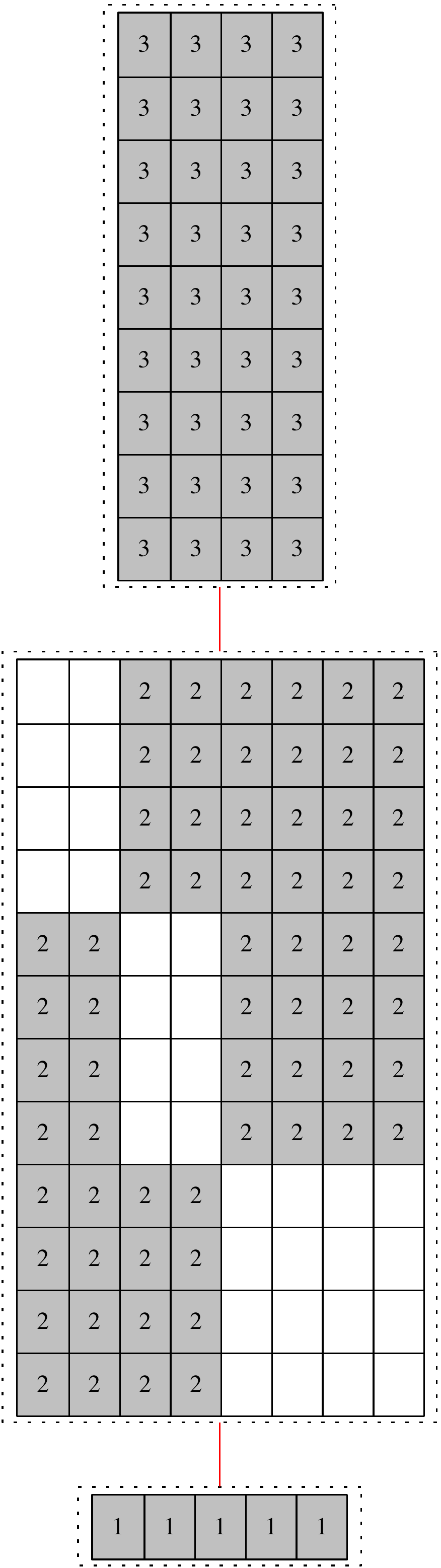}
\ 
\includegraphics[height=92mm]{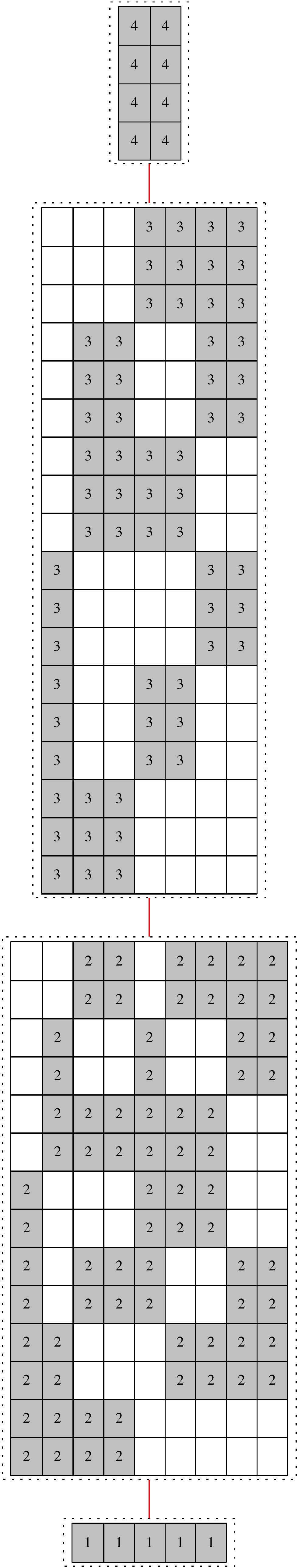}
    \caption{Egg box diagrams of the regular subsemigroups $P=\Reg(\T_5^a)$ in the cases (from left to right): $a=[1,1,1,1,1]$, $a=[1,2,2,2,2]$, $a=[1,1,2,2,2]$, $a=[1,2,3,3,3]$, $a=[1,2,2,3,3]$, $a=[1,2,3,4,4]$.}
    \label{fig:RegTXa}
   \end{center}
\end{figure}

Figures \ref{fig:V4_1233} and \ref{fig:V4_1122_1222} picture the variant $\T_4^a$ with respect to various transformations $a\in\T_4$, and one may see the regular subsemigroup $\Reg(\T_4^a)$ in each case as the collection of $\gDa$-classes containing groups (shaded cells).  Figure \ref{fig:RegTXa} pictures $\Reg(\T_5^a)$ for various choices of $a\in\T_5$ with $\rank(a)\leq4$.
%
%
When one compares Figure~\ref{fig:RegTXa} with Figure \ref{fig:T1...T4}, which pictures the semigroups $\T_1$, $\T_2$, $\T_3$, $\T_4$, a striking pattern seems to emerge.  In each case, $P=\Reg(\T_5^a)$ looks like some kind of ``inflation'' of $\T_r$ (where $r=\rank(a)$), in the sense that one may begin with an egg box diagram of $\T_r$ and then subdivide the cells in some way to obtain an egg box diagram of~$P$; further, it appears that the subdivision is done in such a way that group (resp., non-group) $\gH$-classes of $\T_m$ become rectangular arrays of group (resp., non-group) $\gHa$-classes of $P$, although the reason for the exact number of subdivisions applied to each cell may not yet be apparent.  One of the goals of this section is to explain the reason for this phenomenon.

Now, Theorem \ref{green_thm} enables us to immediately describe Green's relations on $P=\RegTXa$.  Since~$P$ is a regular subsemigroup of $\TXa$, the $\gR$, $\gL$, $\gH$ relations on $P$ are just the restrictions of the corresponding relations on $\TXa$ (see for example \cite{Hig,Howie}), and it is easy to check that this is also true for the ${\gD}={\gJ}$ relation in this case.  So if $\gK$ is one of $\gR$, $\gL$, $\gH$, $\gD$, we will continue to write $\gKa$ for the $\gK$ relation on~$P$, and write $K_f^a$ for the $\gKa$-class of $f$ in $P$ for any $f\in P$.  

\ms
\begin{cor}\label{cor_green}
If $f\in P$, then 
\bit
\itemit{i} $R_f^a =R_f\cap P = \set{g\in P}{\ker(f)=\ker(g)}$, 
\itemit{ii} $L_f^a  =L_f\cap P = \set{g\in P}{\im(f)=\im(g)}$, 
\itemit{iii} $H_f^a  =H_f\cap P = \set{g\in P}{\ker(f)=\ker(g)\text{\emph{ and }}\im(f)=\im(g)}$, 
\itemit{iv} $D_f^a  =D_f\cap P = \set{g\in P}{\rank(f)=\rank(g)}$.
\eit
The $\gDa$-classes of $P$ form a chain: $D_1^a<\cdots<D_r^a$, where $D_m^a=\set{f\in P}{\rank(f)=m}$ for each $m\in\br$. \epfres
\end{cor}

Corollary \ref{cor_green} gives a descriptive characterisation of Green's relations on $P=\RegTXa$; in particular, it relates each relation $\gKa$ on $P$ directly to the relation $\gK$ on $\T_X$.  
But it says nothing about why $P$ appears to be an inflated version of $\T_r$.  In order to explain this phenomenon, we must further explore the structure of~$P$.  We will do this by examining a certain relationship between $P$ and $\T_A$, the full transformation semigroup on $A=\im(a)$, as well as some other well-known subsemigroups of $\T_X$.  Recall that the sets
\[
\TXA=\set{f\in\T_X}{\im(f)\sub A} \AND \TXal=\set{f\in\T_X}{\ker(f)\supseteq\al}
\]
are subsemigroups of $\T_X$.  These semigroups have been studied extensively in the literature, where they are typically referred to as semigroups of transformations of restricted range or restricted kernel (respectively); see for example \cite{Sanwong2011,SS2008,SS2013,MGS2010,Symons1975}, and references therein.  

\ms
\begin{rem}\label{rem_TXA_TXal}
Note that $\TXA=\T_Xa$ and $\TXal=a\T_X$, as subsemigroups of $\T_X$ (with respect to the usual operation).  Indeed, the maps
\[
\rho_a:\TXa\to\TXA=\T_Xa:f\mt fa \AND \lam_a:\TXa\to\TXal=a\T_X:f\mt af
\]
are easily seen to be epimorphisms.  Since products in $\TXA=\T_Xa$ and $\TXal=a\T_X$ are found by forming expressions such as $faga$ and $afag$ (respectively), it should be no surprise that these semigroups play a role in an investigation of the structure of $\TXa$.
Since we are assuming $a$ is an idempotent, it also follows that $\TXA=\T_X\star a$ and $\TXal=a\star \T_X$, as subsemigroups of $\TXa$ (with respect to the $\star$ operation).  As noted in \cite{GMbook}, if $S$ is either $\TXA$ or $\TXal$, the semigroups $S^a$ and $S$ are precisely the same object; that is, $f\star g=fg$ for all $f,g\in S$.  (This is because $a$, being an idempotent of $\T_X$, is a mididentity of both $a\T_X$ and $\T_Xa$.)
\end{rem}

The regular elements of the semigroups $\TXA$ and $\TXal$ have been described in \cite{SS2008} and \cite{MGS2010}, respectively; in terms of our notation, the description is as follows.  Recall that
\[
P_1=\set{f\in\T_X}{\text{$\al$ separates $\im(f)$}} \AND P_2=\set{f\in\T_X}{\text{$A$ saturates $\ker(f)$}}.
\]

\ms
\begin{prop}[Sanwong and Sommanee \cite{SS2008}; Mendes-Gon\c{c}alves and Sullivan \cite{MGS2010}]\label{prop_RegTXA_TXal} \leavevmode\newline
The regular elements of $\TXA$ and $\TXal$ are precisely the sets
\[
\epfreseq
\RegTXA = \TXA\cap P_2 \AND \RegTXal = \TXal\cap P_1.
\]
\end{prop}

The next two propositions are the main structural results of this section.

\ms
\begin{prop}\label{mono_prop}
There is a well-defined monomorphism
\[
\psi : \RegTXa\to \RegTXA\times\RegTXal : f\mt (fa,af).
\]
The image of $\psi$ is the set
\[
\im(\psi) = \bigset{(g,h)\in\RegTXA\times\RegTXal}{\rank(g)=\rank(h),\ g|_A=(ha)|_A}.
\]
In particular, $\RegTXa$ is (isomorphic to) a subdirect product of $\RegTXA$ and $\RegTXal$.
\end{prop}

\pf Let $f\in P=\RegTXa$.  Since $P$ is a subsemigroup of $\TXa$, we have $fa=faa=f\star a\in P$; in particular, $fa\in \T_Xa\cap P_2=\RegTXA$.  A similar calculation shows that $af\in\RegTXal$.  If $f,g\in P$, then
$
(f\star g)\psi = (fag)\psi = ((fag)a,a(fag)) = (fa,af)(ga,ag)=(f\psi)(g\psi),
$
so $\psi$ is a homomorphism.  Suppose now that $f,g\in P$ are such that $f\psi=g\psi$.  So $fa=ga$ and $af=ag$, and we must show that $f=g$.  Since~$A$ saturates $\ker(f)$ and $\ker(g)$, it suffices to show that $\ker(f)=\ker(g)$ and $f|_A=g|_A$.  Now, for any $x\in A$, we have $xf=xaf=xag=xg$, so $f|_A=g|_A$.  Also note that since $f\in P_1$, $\ker(fa)=\ker(f)$.  Similarly, $\ker(ga)=\ker(g)$.  Since $fa=ga$, it follows that $\ker(f)=\ker(g)$.  As noted above, this completes the proof that $\psi$ is injective.

To prove the statement concerning $\im(\psi)$, first suppose $f\in P$ and put $g=fa$ and $h=af$.  Since $f\in P=P_1\cap P_2$, Proposition \ref{reg_prop} gives $\rank(g)=\rank(f)=\rank(h)$.  Since $a$ maps $A$ identically, it follows that $(aq)|_A=q|_A$ for all $q\in\T_X$.  In particular, $(ha)|_A=(afa)|_A=(fa)|_A=g|_A$.
Conversely, suppose $g\in\RegTXA$ and $h\in\RegTXal$ satisfy $\rank(g)=\rank(h)$ and $g|_A=(ha)|_A$.  Put $m=\rank(g)$, and write $g = \trans{G_1 & \cdots & G_m}{a_{k_1} & \cdots & a_{k_m}}$ and $h = \trans{H_1 & \cdots & H_m}{h_1 & \cdots & h_m}$,
noting that $\im(g)\sub A$.  Also, since $h\in\TXal$, there is a partition $\br=I_1\sqcup\cdots\sqcup I_m$ such that $H_j=\bigcup_{i\in I_j}A_i$ for each $j$.  Now, since $g\in\RegTXA$, $A$ saturates $\ker(g)$, so it follows that $G_i\cap A\not=\emptyset$ for all $i$.  Thus,
$g|_A =  \trans{G_1\cap A & \cdots & G_m\cap A}{a_{k_1} & \cdots & a_{k_m}}$.
For each $i\in\bm$, let $l_i\in\br$ be such that $h_i\in A_{l_i}$.  Since $h\in\RegTXal$, $\al$ separates $\im(h)$, so $l_1,\ldots,l_m$ are distinct.  It follows that
$
ha = \trans{H_1 & \cdots & H_m}{a_{l_1} & \cdots & a_{l_m}}$.
Since each $H_i\cap A$ is non-empty (as $H_i$ is a union of $\al$-classes, each of which contains an element of $A$), we have
$(ha)|_A = \trans{H_1\cap A & \cdots & H_m\cap A}{a_{l_1} & \cdots & a_{l_m}}$.
But $g|_A=(ha)|_A$, so (reordering if necessary), it follows that $l_i=k_i$ and $H_i\cap A=G_i\cap A$ for each $i$.  In particular, $h_i\in A_{l_i}=A_{k_i}$ for each $i$.  Now put
$
f = \trans{G_1 & \cdots & G_m}{h_1 & \cdots & h_m}$.
Since $\ker(f)=\ker(g)$ and $\im(f)=\im(h)$, we see that $f\in P$.  It is clear that $fa=g$.  We also have $af=h$ since, for all $j$,
\[
H_jaf=\Big(\bigcup_{i\in I_j}A_i\Big)af=\set{a_i}{i\in I_j}f=(H_j\cap A)f=(G_j\cap A)f=h_j.
\]
It follows that $(g,h)=f\psi$.
Finally, suppose $g\in\RegTXA$ and $h\in\RegTXal$.  To prove the statement about $\RegTXa$ being a subdirect product, we must show that there exist $h'\in\RegTXal$ and $g'\in\RegTXA$ such that $(g,h'),(g',h)\in\im(\psi)$.  First note that $g\in P_2$ by Proposition~\ref{prop_RegTXA_TXal}.  But also $\TXA\sub P_1$, so $g\in P$, and $(g,ag)=(ga,ag)=g\psi$, so we may take $h'=ag$.  Similarly, $h\in P$ and $(ha,h)=h\psi$, and we take $g'=ha$.  This completes the proof. \epf

The homomorphism $\psi$ from the previous result is built up out of the two homomorphisms
\[
\psi_1:P\to \RegTXA:f\mt fa \AND \psi_2:P\to \RegTXal:f\mt af,
\]
which are the restrictions to $P=\RegTXa$ of the epimorphisms $\lam_a$ and $\rho_a$ from Remark \ref{rem_TXA_TXal}.  The last paragraph of the previous proof shows that $\psi_1$ and $\psi_2$ are epimorphisms, indeed projections, since $P$ contains both $\RegTXA$ and $\RegTXal$, and $\psi_1$ (resp., $\psi_2$) maps $\RegTXA$ (resp., $\RegTXal$) identically.

\ms
\begin{prop}\label{epi_prop}
The maps
\[
\phi_1:\RegTXA\to\T_A:g\mt g|_A \AND \phi_2:\RegTXal\to\T_A:g\mt (ga)|_A
\]
are epimorphisms, and the following diagram commutes:
\[
\includegraphics{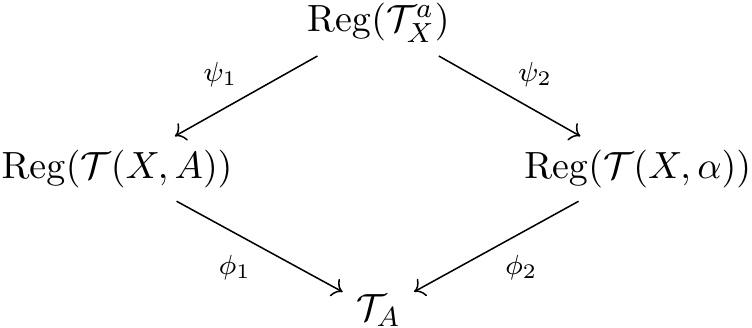}
\]
Further, the induced map $\RegTXa\to\T_A$ is an epimorphism.
\end{prop}

\pf Clearly $\phi_1$ and $\phi_2$ map their domains into $\T_A$.  Note that for any $f\in\T_X$, $f|_A=\id_A\circ f$, where $\id_A$ is the restriction of the identity map to $A$, and $\circ$ denotes the usual composition of partial functions.  So, if $g,h\in\RegTXA$, then 
\[
(gh)\phi_1 = \id_A\circ g\circ h = \id_A\circ g\circ\id_A\circ h = (g\phi_1)(h\phi_1),
\]
since $g=g\circ\id_A$ as $\im(g)\sub A$.  If $g,h\in\RegTXal$, then
\[
(gh)\phi_2 = \id_A\circ g\circ h\circ a = \id_A\circ g\circ a\circ h\circ a = \id_A\circ g\circ a\circ\id_A\circ h\circ a = (g\phi_2)(h\phi_2),
\]
since $h=a\circ h$ as $\al\sub\ker(h)$, and $a=a\circ\id_A$.  So $\phi_1$ and $\phi_2$ are homomorphisms.  That the diagram commutes follows from the fact that $(afa)|_A=(fa)|_A$ for all $f\in P$, as observed in the proof of Proposition~\ref{mono_prop}.  Finally, let $q\in\T_A$, and write
$q =  \trans{a_1 & \cdots & a_r}{a_{k_1} & \cdots & a_{k_r}}$.
(This notation is not supposed to imply that $k_1,\ldots,k_r$ are distinct.)  Put
$f =  \trans{A_1 & \cdots & A_r}{a_{k_1} & \cdots & a_{k_r}}$.
Then clearly, $f\in\RegTXA\cap\RegTXal$ and $q=f\phi_1=f\phi_2$, showing that $\phi_1$ and $\phi_2$ are surjective.  Note that, in fact, $f\in P$ and $q=f\phi_1=(fa)\phi_1=f(\psi_1\phi_1)$, showing that $\psi_1\phi_1$ is surjective, and completing the proof. \epf

\begin{rem}
The previous result displays the structure of $P=\RegTXa\cong\im(\psi)$ as a \emph{pullback product} of $\RegTXA$ and $\RegTXal$ with respect to $\T_A$.  Namely, $\im(\psi)$ consists of all pairs $(g,h)$ such that $g\phi_1=h\phi_2$.  Pullback products have been studied in various contexts in universal algebra and semigroup theory (where they are sometimes referred to as \emph{spined products}); see for example \cite{Fleischer1955,CB1998,CPB2008,Fuchs1952,Kimura1958}.
\end{rem}

From now on, we will denote by $\phi=\psi_1\phi_1=\psi_2\phi_2$ the epimorphism $P\to\T_A:f\mt(fa)|_A$.
%
If $f\in P$, we will write $\fb=f\phi\in\T_A$.  If $U\sub P$, we write $\Ub=\set{\ub}{u\in U}\sub\T_A$.

We now show how $\phi:P\to\T_A$ may be used to relate Greens relations on the semigroups $P$ and $\T_A$.
If $f,g\in P$ and $\gK$ is one of $\gL$, $\gR$, $\gH$, $\gD$, we say $f\gKh g$ if $\fb\gK\gb$ in $\T_A$.  Denote by $\Kh_f=K_{\fb}\phi^{-1}$ the $\gKh$-class of $f$ in $P$.  
Recall that $\lam_i=|A_i|$ for each $i\in\br$, and that $\Lam_I=\lam_{i_1}\cdots\lam_{i_m}$ if $I=\{i_1,\ldots,i_m\}\sub\br$.  If $Y$ is a set and $0\leq m\leq |Y|$, we write ${Y\choose m}$ for the set of all $m$-element subsets of $Y$.  Recall that a \emph{rectangular band} is a semigroup of the form $I\times J$ with product $(i_1,j_1)(i_2,j_2)=(i_1,j_2)$, and that a \emph{rectangular group} is a direct product of a rectangular band with a group.

\ms
\begin{thm}\label{inflation_thm}
Let
$
f=\trans{F_1 & \cdots & F_m}{f_1 & \cdots & f_m}\in P$,
where $m=\rank(f)$ and $f_i\in A_{k_i}$ for each $i$, and put $I=\{k_1,\ldots,k_m\}$.
\bit
\itemit{i} $\Rh_f$ is the union of $m^{n-r}$\ $\gRa$-classes of $P$.
\itemit{ii} $\Lh_f$ is the union of $\Lam_I$ $\gLa$-classes of $P$.
\itemit{iii} $\Hh_f$ is the union of $m^{n-r}\Lam_I$ $\gHa$-classes of $P$, each of which has size $m!$.  The map $\phi:P\to\T_A$ is injective when restricted to any $\gHa$-class of $P$.
\itemit{iv} If $H_{\fb}$ is a non-group $\gH$-class of $\T_A$, then each $\gHa$-class of $P$ contained in $\Hb_f$ is a non-group.
\itemit{v} If $H_{\fb}$ is a group $\gH$-class of $\T_A$, then each $\gHa$-class of $P$ contained in $\Hb_f$ is a group isomorphic to the symmetric group $\S_m$.  Further, $\Hb_f$ is a rectangular group; specifically, $\Hb_f$ is isomorphic to a direct product of an $m^{n-r}\times\Lam_I$ rectangular band with $\S_m$.
\itemit{vi} ${\gDh}={\gDa}$, so $\Dh_f=D_f^a=D_m^a=\set{g\in P}{\rank(g)=m}$ is the union of:
\begin{itemize}
\itemit{a} $m^{n-r}S(r,m)$ $\gRa$-classes of $P$,
\itemit{b} $\sum_{J\in{\br \choose m}}\Lam_J$ $\gLa$-classes of $P$,
\itemit{c} $m^{n-r}S(r,m)\sum_{J\in{\br \choose m}}\Lam_J$ $\gHa$-classes of $P$.
\eit
\eit
\end{thm}

\pf First observe that if $\rho:S\to T$ is an epimorphism of semigroups, and if $K$ is a $\gK$-class of $T$ where $\gK$ is one of $\gL,\gR,\gH$, then $K\rho^{-1}$ is a union of $\gK$-classes of $S$.
\bit
\item[(i)] By the above observation, it suffices to count the number of $\gRa$-classes contained in $\Rh_f$.  An $\gRa$-class $R_g^a$ contained in $\Rh_f$ is completely determined by the common kernel of all its members, namely $\ker(g)$.  Such a kernel is constrained so that it has $m$ equivalence classes and $\ker(\gb)=\ker(\fb)=(F_1\cap A|\cdots|F_m\cap A)$.  To construct $\ker(g)$ from $\ker(\gb)$, the remaining $n-r$ elements of $X\sm A$ may be assigned to the $m$ $\ker(\gb)$-classes arbitrarily, and there are $m^{n-r}$ ways to do this.

\item[(ii)] An $\gLa$-class $L_g^a$ contained in $\Lh_f$ is completely determined by the common image of all its members, namely $\im(g)$.  Such an image is constrained so that it has size $m$ and $\im(\gb)=\im(\fb)=\{a_{k_1},\ldots,a_{k_m}\}$.  So $\im(g)$ must contain one element of $A_{k_i}$ for each $i$, and may be chosen in $\lam_{k_1}\cdots\lam_{k_m}=\Lam_I$ ways.

\item[(iii)] The statement concerning the number of $\gHa$-classes contained in $\Hb_f$ follows immediately from (i) and (ii).  By Theorem \ref{green_thm}, $H_f^a=H_f$, so $|H_f^a|=m!$.  If $g\gHa f$, then 
$
g=\trans{F_1 & \cdots & F_m}{f_{1q} & \cdots & f_{mq}}
$ and $
\gb=\trans{F_1\cap A & \cdots & F_m\cap A}{a_{k_{1q}} & \cdots & a_{k_{mq}}}
$
for some $q\in\S_m$.  So it follows that $\phi$ is injective when restricted to $H_f^a$.  Since $H_f^a$ is an arbitrary $\gHa$-class of $\RegTXa$, the proof of (iii) is complete.

\item[(iv)] If $H_{\fb}$ is a non-group $\gH$-class, then $\gb^2\not\in H_{\fb}$ for all $g\in\Hh_f$.  Since $\im(\gb^2)\sub\im(\gb)$ and $\ker(\gb^2)\supseteq\ker(\gb)$, it then follows that $\rank(g^2)=\rank(\gb^2)<\rank(\gb)=\rank(g)$, so $g^2\not\in H_g^a$, whence $H_g^a$ is a non-group $\gHa$-class of $P$.  


\item[(v)] Suppose $H_{\fb}$ is a group.  Then $\gb^2\in H_{\fb}$ for any $g\in \Hh_f$, so $\rank(g^2)=\rank(\gb^2)=\rank(\gb)=\rank(g)$.  But $\im(g^2)\sub\im(g)$ and $\ker(g^2)\supseteq\ker(g)$, so it follows that $\im(g^2)=\im(g)$ and $\ker(g^2)=\ker(g)$, whence $g^2\gHa g$, whence $H_g^a$ is a group.
By (iii), the restriction of $\phi$ to $H_f^a$ yields an isomorphism onto $H_{\fb}\cong\S_m$.

Consider an arbitrary element $g\in\Hh_f$, and write $\ker(g)=(G_1|\cdots|G_m)$ and $\im(g)=\{g_1,\ldots,g_m\}$ where $a_{k_i}\in G_i$ and $g_i\in A_{k_i}$ for each $i\in\bm$.  Then there is a permutation $p_g\in\S_m$ such that
$
g=\trans{G_1 & \cdots & G_m}{g_{1p_g} & \cdots & g_{mp_g}}$.
In this way, we see that $g$ is completely determined by $\ker(g)$, $\im(g)$ and $p_g$, and we write $g\equiv[\ker(g),\im(g),p_g]$.   If 
$
h=\trans{H_1 & \cdots & H_m}{h_{1p_h} & \cdots & h_{mp_h}}
\equiv[\ker(h),\im(h),p_h]
$
is another element of~$\Hh_f$, then we have $g\star h=[\ker(g),\im(h),p_gp_h]$.  Indeed, we have $\ker(g\star h)=\ker(g)$ and $\im(g\star h)=\im(h)$, as $g\star h\in\Hh_f$ gives $\rank(g\star h)=m$, and if $x\in G_i$ is arbitrary, then
\[
x \lmap{g} g_{ip_g} \lmap{a} a_{k_{ip_g}} \lmap{h} h_{ip_gp_h}.
\]
Now let $K=\set{\ker(g)}{g\in\Hh_f}$ and $I=\set{\im(g)}{g\in\Hh_f}$.  Then $K\times I$ is a rectangular band under the product
$
(\be,B)(\ga,C)=(\be,C),
$
and by the above rule for multiplication in $\Hh_f$, we immediately see that the map
\[
\Hh_f\to K\times I\times \S_m: g\equiv[\ker(g),\im(g),p_g]\mt(\ker(g),\im(g),p_g)
\]
is an isomorphism. The dimensions of the rectangular band are given by parts (i) and (ii), above. 

\item[(vi)] We deduce ${\gDh}={\gDa}$ immediately from the fact that $\rank(\fb)=\rank(f)$ for all $f\in P$. 
%
The number of $\gRh$-classes in $D_m^a$ is equal to the number of $\gR$-classes in $D_m\sub\T_A$, which is equal to $S(r,m)$; (a) now follows from (i).  
Part (b) follows from (ii) and the fact that the $\gL$-classes contained in $D_m\sub\T_A$ (and hence the $\gLh$-classes contained in $D_m^a$) are indexed by the $m$-element subsets of $A$.  Part (c) follows immediately from (a) and (b).
%
%
\epf
\eit

\begin{rem}
See also \cite[Proposition 3.1]{Tsyaputa2004} for formulae for the number of singleton $\gRa$- and $\gLa$-classes of~$\TXa$, and various other parameters.
\end{rem}

As an immediate consequence of Theorem \ref{inflation_thm}, we may give the size of $P=\RegTXa$.  

\ms
\begin{cor}\label{cor_|P|}
We have
$\displaystyle{
|\RegTXa| = \sum_{m=1}^r m!m^{n-r}S(r,m) \sum_{I\in{\br \choose m}} \Lam_I. 
}$
\end{cor}

\pf From parts (vi) and (iii) of Theorem \ref{inflation_thm}, we have $|D_m^a|=m!m^{n-r}S(r,m) \sum_{I\in{\br \choose m}} \Lam_I$ for each $m\in\br$.  Summing over all $m$ gives the result. \epf

The top $\gDa$-class of $P$ is the set
\[
D_r^a = \S_A\phi^{-1} = \set{f\in P}{\rank(f)=r}.
\]
We will write $D=D_r^a$ for this set.  As a special case of Theorem \ref{inflation_thm}(v), $D$ is a rectangular group; it is isomorphic to the direct product of an $r^{n-r}\times\Lam$ rectangular band with the symmetric group $\S_r$.  (Recall that $\Lam=\lam_1\cdots\lam_r$.)  Since $D$ is the pre-image of $\S_A$ under the map $\phi:P\to\T_A$, we may think of $D$ as a kind of ``inflation'' of $\S_A$, the group of units of $\T_A$.  In fact, we will soon see that $D=\RP(P)$ is precisely the set of regularity preserving elements of $P$, so that $D$ may be thought of as an alternative to the group of units in the non-monoid $P$, as noted in Section \ref{sect:variants}.


In order to avoid confusion when discussing idempotents, if $U\sub\T_X$, we will write
\[
E(U)=\set{f\in U}{f=f\circ f} \AND \Ea(U)=\set{f\in U}{f=f\star f}
\]
for the set of idempotents from $U$ with respect to the different operations on $\T_X$ and $\TXa$.  
%
%
%
Recall that an element $u$ of a semigroup $S$ is a mididentity if $xuy=xy$ for all $x,y\in S$.

\ms
\begin{lemma}\label{lemma_aea}
Let $e\in \Ea(D)$.  Then $aea=a$.  In particular, $e$ is a mididentity for both $\TXa$ and $P$.
\end{lemma}

\pf Since $\rank(e)=r$, we may write
$
e=\trans{E_1 & \cdots & E_r}{e_1 & \cdots & e_r}$.
Since $\al$ separates $\im(e)=\{e_1,\ldots,e_r\}$, we may assume (reordering if necessary) that $e_i\in A_i$ for each $i$.  It follows that
$
ea=\trans{E_1 & \cdots & E_r}{a_1 & \cdots & a_r}$.
Since $e\in \Ea(D)$, we see that $e=e\star e=eae$.  It follows that $ea=eaea$, so $ea\in E(\T_X)$, whence $a_i\in E_i$ for each $i$.  It follows then that
$
aea=\trans{A_1 & \cdots & A_r}{a_1 & \cdots & a_r} = a$.  If $f,g\in\T_X$, then $f\star e\star g=faeag=fag=f\star g$, showing that $e$ is a mididentity for $\TXa$ (and hence also for $P\sub\TXa$) and
completing the proof. \epf

\begin{prop}
The top $\gDa$-class, $D=D_r^a$, of $P=\RegTXa$ is precisely the set $\RP(P)$ of all regularity preserving elements of $P$.
\end{prop}

\pf By Proposition \ref{RPS_prop}(i), it suffices to show that $\Ea(\RP(P))=\Ea(D)$.  By Proposition \ref{RPS_prop}(iii) and Lemma~\ref{lemma_aea}, we see that $\Ea(D)\sub\Ea(\RP(P))$.  Conversely, suppose $e\in\Ea(P)\sm\Ea(D)$.  Then $\rank(e)<r$, and so if $f\in\Ea(D)$ is arbitrary, then $\rank(f\star e)=\rank(fae)\leq\rank(e)<r$, so $f\star e$ does not belong to $D=D_f^a$ and, in particular, $f\star e$ is not $\gRa$-related to $f$, from which we deduce from Proposition \ref{RPS_prop}(ii) that $e\not\in\RP(P)$.  This shows that $\Ea(\RP(P))\sub \Ea(D)$, and completes the proof. \epf



Our next goal is to calculate the rank of $P=\RegTXa$.  Recall that the \emph{relative rank}, denoted $\rank(S:U)$, of a semigroup $S$ with respect to a subset $U\sub S$ is defined to be the minimum cardinality of a subset $V\sub S$ such that $S=\la U\cup V\ra$.  The concept of relative rank was first introduced in \cite{HHR}, and has played a major role in a number of investigations \cite{HHR2,HHMR,HMR,HMMR,CMMP,AM07}. 




\ms
\begin{lemma}\label{lemma_rankP_D_PD}
We have $\rank(P)=\rank(D)+\rank(P:D)$.
\end{lemma}

\pf This follows quickly from the fact that $D$ is a subsemigroup of $P$ and $P\sm D$ an ideal.~\epf


The next result may be easily be proved directly, but it is a special case of \cite[Theorem 4.7]{Ruskuc1994} (see also \cite{Gray2005}) so we omit the proof.

\ms
\begin{lemma}[Ru\v skuc \cite{Ruskuc1994}]\label{lemma_IJG}
Let $I$ and $J$ be non-empty sets, and $G$ a group.  Let $S=I\times J\times G$ be the rectangular group with product defined by
$
(i_1,j_1,g_1)(i_2,j_2,g_2)=(i_1,j_2,g_1g_2).
$
Then
\[\epfreseq
\rank(S)=\max\big\{|I|,|J|,\rank(G)\big\}.
\]
\end{lemma}



We wish to apply Lemma \ref{lemma_IJG} to calculate the rank of the rectangular group $D$.  To do this, we need to calculate $\max\{r^{n-r},\Lam\}$.  Recall that we are assuming $1<r<n$.  

\ms
\begin{lemma}\label{lemma_r^n-r}
We have $r^{n-r}\geq\Lam=\lam_1\cdots\lam_r$.
\end{lemma}

\pf First note that if $r=2$ and $n=3$, then we must have $\{\lam_1,\lam_2\}=\{1,2\}$, in which case $r^{n-r}=\Lam=2$.  Now suppose $(r,n)\not=(2,3)$.  Elementary calculus shows that the maximum value of the product $x_1\cdots x_r$, where $x_1+\cdots+x_r=n$ and $x_1,\ldots,x_r\geq0$ are real numbers, occurs when $x_1=\cdots=x_r=n/r$.  It follows that $\Lam\leq (n/r)^r=n^r/r^r$.  So it suffices to show that $n^r/r^r\leq r^{n-r}=r^n/r^r$, which is equivalent to $n^r\leq r^n$.  This, in turn, is equivalent to $\ln(n)/n\leq \ln(r)/r$.  Now, $f(x)=\ln(x)/x$ is a decreasing function for $x>e\approx2.718$.  In particular, $f(3)>f(4)>f(5)>\cdots$, so the result holds for $r\geq3$.  We also have $f(2)=f(4)$, so the result holds for $r=2$ and $n\geq4$.  We have already covered the case $(r,n)=(2,3)$. \epf

\begin{cor}\label{cor_rankD}
We have $\rank(D)=r^{n-r}$.
\end{cor}

\pf Recall that $D$ is isomorphic to the direct product of an $r^{n-r}\times\Lam$ rectangular band with the symmetric group $\S_r$.  So Lemma \ref{lemma_IJG} gives
$
\rank(D) = \max\big\{ r^{n-r}, \Lam, \rank(\S_r) \big\}.
$
We have already seen that $r^{n-r}\geq\Lam$.  Also, $\rank(\S_2)=1$, while $\rank(\S_r)=2$ if $r\geq3$.  So it follows that $r^{n-r}\geq\rank(\S_r)$.~\epf

The next technical lemma will help us calculate $\rank(P:D)$.  It is quite a bit stronger than we need at this point (we only require the $m=r$ case at the moment), but we will use the full strength in subsequent sections when we consider ideals and the idempotent generated subsemigroup of $\TXa$.

\ms
\begin{lemma}\label{lemma_e1e2}
Suppose $f,g\in P$ are such that $\fb=\gb$.  Then for any $\rank(f)\leq m\leq r$, there exist idempotents $e_1,e_2\in \Ea(D_m^a)$ such that $f=e_1\star g\star e_2$.
\end{lemma}

\pf Put $l=\rank(f)=\rank(g)$ and write
$
f = \trans{F_1 & \cdots & F_l}{f_1 & \cdots & f_l}
$ and $
g = \trans{G_1 & \cdots & G_l}{g_1 & \cdots & g_l}$,
where $f_i\in A_{k_i}$ for each~$i$.  Since
$
\trans{F_1\cap A & \cdots & F_l\cap A}{a_{k_1} & \cdots & a_{k_l}}
=\fb=\gb$,
we may assume (reordering if necessary) that $g_i\in A_{k_i}$ for all $i$, in which case also $G_i\cap A=F_i\cap A$.  Let $\br\sm\{k_1,\ldots,k_l\}=\{j_1,\ldots,j_{r-l}\}$, and put $B=A_{j_{m-l+1}}\cup\cdots\cup A_{j_{r-l}}$.  (Note that $B=\emptyset$ if $m=r$.)  Define
\[
e_2 = \bigtrans{A_{k_1}\cup B & A_{k_2} & \cdots & A_{k_l} & A_{j_1} & \cdots & A_{j_{m-l}}}{f_1 & f_2 & \cdots & f_l & a_{j_1} & \cdots & a_{j_{m-l}}}.
\]
For each $s\in\bl$, let $F_s\cap A=\{a_{i_{s1}},\ldots,a_{i_{sq_s}}\}$, noting that $F_s\cap A\not=\emptyset$ and $q_1+\cdots+q_l=r$.  For each $s$, choose $1\leq p_s\leq q_s$ such that $p_1+\cdots+p_l=m$, and choose a partition $F_s=F_{s1}\sqcup\cdots\sqcup F_{sp_s}$ so that $a_{i_{st}}\in F_{st}$ for each $t$.  
Define
\[
e_1 = \bigtrans{F_{11} & \cdots & F_{1p_1} & \cdots & F_{l1} & \cdots & F_{lp_l}}{a_{i_{11}} & \cdots & a_{i_{1p_1}} & \cdots & a_{i_{l1}} & \cdots & a_{i_{lp_l}}}.
\]
One may easily check that $e_1,e_2\in E(\T_X)$.  Since also $e_1a=e_1$ and $ae_2=e_2$, it follows that $e_1,e_2\in\Ea(D_m^a)$.  Now let $x\in F_s$ be arbitrary.  Then $xe_1\in F_s\cap A=G_s\cap A$, so 
\[
x  \lmap{e_1} xe_1 \lmap{a} xe_1 \lmap{g} g_s \lmap{a} a_{k_s} \lmap{e_2} f_s=xf,
\]
showing that $f=e_1\star g\star e_2$, as desired. \epf

\begin{lemma}\label{lemma_rankPD}
If $f\in D_{r-1}^a$ is arbitrary, then $P=\la D\cup\{f\}\raa$.  Consequently, $\rank(P:D)=1$.
\end{lemma}

\pf Since $\la D\raa=D\not=P$ (as $r>1$), $\rank(P:D)\geq1$ so it suffices to prove the first statement.
Note that $\Db=\set{\gb}{g\in D}$ is equal to $\S_A$, and $\fb\in\T_A$ satisfies $\rank(\fb)=r-1$.  It follows that $\T_A=\la\Db\cup\{\fb\}\ra$.  Now let $g\in P$ be arbitrary.  Choose $h_1,\ldots,h_k\in D\cup\{f\}$ such that $\gb=\hb_1\cdots\hb_k$, and put $h=h_1\star\cdots\star h_k\in \la D\cup\{f\}\raa$.  Then $\hb=\gb$, so Lemma \ref{lemma_e1e2} tells us that $g=e_1\star h\star e_2\in \la D\cup\{f\}\raa$ for some $e_1,e_2\in\Ea(D)$. \epf

As an immediate consequence of Lemmas \ref{lemma_rankP_D_PD} and \ref{lemma_rankPD} and Corollary \ref{cor_rankD}, we have the following.

\ms
\begin{thm}\label{thm_rankP}
If $1<r<n$, then $\rank(\RegTXa)=r^{n-r}+1$. \epf
\end{thm}

\begin{rem}\label{rem_Rank_RegTXA}
It was shown in \cite[Theorem 3.6]{SS2013} that $\rank(\RegTXA)=r^{n-r}+1$, also.  See also Theorem~\ref{thm_RegTXa_ideals} and Remark \ref{rem_Rank_RegTXA2}.
%
If $r=1$, then $\RegTXa=D_1$ is an $n$-element right zero semigroup, so we have $\rank(\RegTXa)=n$ in this case.  If $r=n$, then $\RegTXa=\TXa\cong\T_X$, so $\rank(\RegTXa)=\rank(\T_X)$, which is equal to $1$ (if $n\leq1$), $2$ (if $n=2$) or $3$ (if $n\geq3$).
\end{rem}

\ms
\begin{rem}
The natural task of classifying and enumerating the generating sets of $P$ of the minimal size $r^{n-r}+1$ seems virtually unassailable.  Indeed, by the proof of Lemma \ref{lemma_IJG} (see \cite[Theorem 4.7]{Ruskuc1994}), such a classification would involve classifying and enumerating all generating sets of $\S_r$ of size at most~$r^{n-r}$.
\end{rem}

\section{The idempotent generated subsemigroup $\la \Ea(\TXa)\raa$}\label{sect:ETXa}


In this section, we investigate the idempotent generated subsemigroup $\EXa=\la\Ea(\TXa)\raa$ of $\TXa$.  Our main results include a proof that $\EXa=\Ea(D)\cup(P\sm D)$, a calculation of $\rank(\EXa)=\idrank(\EXa)$, and an enumeration of the idempotent generating sets of this minimal possible size.
Since the solution to every problem we consider is trivial when $r=1$, and well-known when $r=n$, we will continue to assume that $1<r<n$.  To simplify notation, we will write $\EaTXa=\Ea(\TXa)=\Ea(P)$, so $\EXa=\la E\raa$.  
We begin with a simple observation; part (ii) is proved in \cite[Proposition 13.3.2]{GMbook}, where the idempotents were characterised in a different way (we include a short proof for completeness).

\ms\ms
\begin{prop}
\bit
\itemit{i} $\ds{\EaTXa=\Ea(\TXa) = \set{f\in\T_X}{(af)|_{\im(f)}=\id_{\im(f)}}}$\emph{;}
\itemit{ii} $\ds{|E| = \sum_{m=1}^r m^{n-m} \sum_{I\in{\br\choose m}}\Lam_I}$.
\eit
\end{prop}

\pf Part (i) is easily checked.  For part (ii), note that to specify an idempotent $f\in E$, we first choose $m=\rank(f)=\rank(\fb)\in\br$, then $\im(\fb)=\{a_{i_1},\ldots,a_{i_m}\}$, then $\im(f)=\{b_1,\ldots,b_m\}$ where $b_k\in A_{i_k}$ for each $k\in\bm$.  Note that the condition $(af)|_{\im(f)}=\id_{\im(f)}$ simply says that $a_{i_k}f=b_k$ for each $k$.  The remaining $n-m$ points of $X\sm\{a_{i_1},\ldots,a_{i_m}\}$ may be mapped arbitrarily by $f$ to any of the points from $\{b_1,\ldots,b_m\}$. Multiplying the number of choices at each step, and adding as appropriate, gives the desired result. \epf

\begin{lemma}\label{lemma_EDm_EDma}
If $f\in E(\T_A)$, then there exists $e\in E=\Ea(\TXa)$ such that $f=\eb$ and $\rank(e)=\rank(f)$.
\end{lemma}

\pf One easily checks that $e=\trans{A_1&\cdots&A_r}{a_1f&\cdots&a_rf}$ satisfies the desired conditions.~\epf

Recall that $\TASA$ is idempotent generated; see Theorem \ref{thm_TXSX}.

\ms
\begin{lemma}\label{lemma_PDsubE}
Let $V\sub \Ea(P\sm D)$ be an arbitrary set of idempotents such that $\TASA=\la\Vb\ra$.  Then $\la \Ea(D)\cup V\raa$ contains $P\sm D$.
\end{lemma}

\pf Let $f\in P\sm D$ be arbitrary.  Choose $e_1,\ldots,e_k\in V$ so that $\fb=\eb_1\cdots\eb_k$, and put $g=e_1\star\cdots\star e_k\in\la V\raa$.  So $\gb=\fb$, and Lemma \ref{lemma_e1e2} tells us that there exist $e_0,e_{k+1}\in \Ea(D)$ such that $f=e_0\star g\star e_{k+1}\in\la \Ea(D)\cup V\raa$.~\epf

We may now describe the idempotent generated subsemigroup $\EXa=\la \Ea(\TXa)\raa$ of $\TXa$.  

\ms
\begin{thm}\label{thm_EXa}
We have $\EXa=\la \EaTXa\raa=\Ea(D)\cup (P\sm D)$, where $\EaTXa=\Ea(\TXa)=\Ea(P)$ and $D=D_r^a$ is the top $\gDa$-class of $P=\RegTXa$.
\end{thm}

\pf First, $\Ea(D)\sub \EaTXa$, and it follows from Lemma \ref{lemma_PDsubE} that $P\sm D\sub\EXa$.  It remains to show that $\EXa\sub \Ea(D)\cup (P\sm D)$.  So suppose $f\in\EXa$, and consider an expression $f=e_1\star\cdots\star e_k$, where $e_1,\ldots,e_k\in E$.  We must show that $f\in \Ea(D)\cup(P\sm D)$.  If $f\in P\sm D$, we are done, so suppose $f\in D$.  Since $P\sm D$ is an ideal, it follows that $e_1,\ldots,e_k\in D$.  But $D$ is a rectangular group, so $\Ea(D)$ is a rectangular band.  In particular, $f=e_1\star\cdots\star e_k\in \Ea(D)$. \epf

\begin{rem}
Theorem \ref{thm_EXa} is a pleasing analogue of Howie's result \cite{Howie1966}
that $\la E(\T_X)\ra = \{1\}\cup(\TXSX)$, since $\{1\}=E(\S_X)$, where $\S_X$ is the top $\gD$-class of $\T_X$.  Also, $\S_X=G(\T_X)=\RP(\T_X)$ and, while $\TXa$ has no group of units as it is not a monoid, it is still the case that $D=\RP(P)$.
%
\end{rem}


Now that we have described the elements of the semigroup $\EXa$, the next natural task is to calculate its rank and idempotent rank.  
%
%
To do this, we need the first part of the next result;
the second part will be of use when we later enumerate the idempotent generating sets of $\EXa$ of minimal possible size.

\ms
\begin{lemma}\label{lemma_IJ}
Let $I$ and $J$ be non-empty sets, and $S=I\times J$ the rectangular band with product defined by
$
(i_1,j_1)(i_2,j_2)=(i_1,j_2).
$
Then
\[
\rank(S)=\idrank(S)=\max\big\{|I|,|J|\big\}.
\]
If $I$ and $J$ are finite, then the number of (idempotent) generating sets of this smallest possible size is equal to $y!S(x,y)$, where $x=\max\big\{|I|,|J|\big\}$ and $y=\min\big\{|I|,|J|\big\}$.  
\end{lemma}

\pf Note that $S$ is (isomorphic to) a rectangular group with respect to a trivial group, which has rank~$1$, so the statement about $\rank(S)=\idrank(S)$ follows immediately from Lemma \ref{lemma_IJG} (or may easily be proved directly).
Now let $U$ be an arbitrary generating set of $S$ of minimal possible size.  By duality, we may assume that $x=|I|$ and $y=|J|$.  By considering an expression $(i,j)=u_1\cdots u_k$, where $u_1,\ldots,u_k\in U$, we see that for each $i\in I$, $U$ contains $(i,j_i)$ for some $j_i\in J$.  
Since we are assuming that $|U|=x=|I|$, we see that in fact $U=\set{(i,j_i)}{i\in I}$.  A similar consideration shows that $J=\set{j_i}{i\in I}$, so $i\mt j_i$ defines a surjective map $I\to J$.  (In fact, considered as a binary relation, $U$ \emph{is} a surjective map $I\to J$.)  Conversely, any surjective map $I\to J$ determines an idempotent generating set of $S$ of size $x=|I|$.  Since the number of surjective functions from an $x$-set to a $y$-set is $y!S(x,y)$, the result follows.~\epf


Since $\Ea(D)$ is an $r^{n-r}\times\Lam$ rectangular band, the next result follows from Lemmas \ref{lemma_r^n-r} and \ref{lemma_IJ}.

\ms
\begin{cor}\label{cor_rankEaD}
We have $\rank(\Ea(D))=\idrank(\Ea(D)) = r^{n-r}$, and the number of minimal (idempotent) generating sets of $\Ea(D)$ is equal to $\Lam!S(r^{n-r},\Lam).$ \epfres
\end{cor}




\ms
\begin{thm}\label{thm_rankEXa}
We have $\rank(\EXa)=\idrank(\EXa)=r^{n-r}+\rho_r$, where $\rho_2=2$ and $\rho_r={r\choose2}$ if $r\geq3$.
\end{thm}

\nc{\Wb}{\overline{W}}

\pf As in Lemma \ref{lemma_rankP_D_PD}, we have $\rank(\EXa)=\rank(\Ea(D))+\rank(\EXa:\Ea(D))$ so, by Corollary~\ref{cor_rankEaD}, it remains to show that:
\bit
\item[(i)] there exists a set $V\sub E$ of size $\rho_r$ such that $\EXa=\la\Ea(D)\cup V\raa$, and
\item[(ii)] if $W\sub \EXa\sm \Ea(D)=P\sm D$ satisfies $\EXa=\la\Ea(D)\cup W\raa$, then $|W|\geq\rho_r$.
\eit
Let $U\sub E(\T_A)$ be an arbitrary idempotent generating set of $\TASA$ with $|U|=\rho_r$.  By Lemma \ref{lemma_EDm_EDma}, we may choose $V\sub E$ such that $|V|=\rho_r$ and $\Vb=U$.
Since $U$ is a generating set of $\TASA$, Lemma \ref{lemma_PDsubE} and Theorem \ref{thm_EXa} give $\la \Ea(D)\cup V\raa=\EXa$, establishing (i).  

Next, suppose $\EXa=\la \Ea(D)\cup W\raa$, where $W\sub P\sm D$.  We will show that $\Wb$ generates $\TASA$.  Indeed, let $g\in\TASA$ be arbitrary, and choose any $h\in P$ such that $\hb=g$.  Since $\rank(h)=\rank(\hb)=\rank(g)<r$, it follows that $h\in P\sm D\sub\EXa$.  Consider an expression $h=u_1\star\cdots\star u_k$, where $u_1,\ldots,u_k\in \Ea(D)\cup W$.  
Now, $g=\hb=\ub_1\cdots\ub_k$.  If any of the $u_i$ belongs to $\Ea(D)$, then $\ub_i=1$, the identity element of $\T_A$.  So the factor $\ub_i$ is not needed in the product $g=\ub_1\cdots\ub_k$.  After cancelling all such factors, we see that $g$ is a product of elements from $\Wb$.  Since $g\in\TASA$ was arbitrary, we conclude that $\TASA=\la\Wb\ra$.  In particular, $|W|\geq|\Wb|\geq\rank(\TASA)=\rho_r$, giving (ii). \epf


Now that we know the size of a minimal (idempotent) generating set for $\EXa$, our next task is to enumerate the idempotent generating sets of this size.
For $i,j\in\br$ with $i\not=j$, let $e_{ij}\in E(\T_r)$ and $\ep_{ij}\in E(\T_A)$ be the transformations of $\br$ and $A$ (respectively) defined by
\[
ke_{ij} = \begin{cases}
i &\text{if $k=j$}\\
k &\text{if $k\in\br\sm\{j\}$}
\end{cases}
\AND
a_k\ep_{ij} = \begin{cases}
a_i &\text{if $k=j$}\\
a_k &\text{if $k\in\br\sm\{j\}$.}
\end{cases}
\]
Note that $a_k\ep_{ij}=a_{ke_{ij}}$ for all $i,j,k$.
Recall that $\bbT_Y$ denotes the set of all strongly connected tournaments on the vertex set $Y$ with $|Y|\geq3$.  We will write $\bbT_r$ for~$\bbT_{\br}$.  Recall also the convention that $\bbT_2=\bbT_{{\bf 2}}$ consists of the single directed graph on vertex set ${\bf2}=\{1,2\}$ with edges $(1,2)$ and $(2,1)$.  If $j\in\br$ and $\Ga\in\bbT_r$, we write $d_\Ga^+(j)$ for the in-degree of vertex $j$ in $\Ga$.  

\ms
\begin{thm}\label{thm_nrIGS}
The number of idempotent generating sets of $\EXa$ of the minimal possible size $r^{n-r}+\rho_r$ is equal to
\[
\big[(r-1)^{n-r}\Lam\big]^{\rho_r}\Lam!S(r^{n-r},\Lam)\sum_{\Ga\in\bbT_r}\frac{1}{\lam_1^{d_\Ga^+(1)}\cdots\lam_r^{d_\Ga^+(r)}}.
\]
\end{thm}

\pf Let $U$ be an arbitrary minimal idempotent generating set of $\EXa=\Ea(D)\cup(P\sm D)$.  Put $U_1=U\cap \Ea(D)$ and $U_2=U\cap(P\sm D)$.  Since $P\sm D$ is an ideal of $\EXa$, it follows that $U_1$ is a (minimal) idempotent generating set of $\Ea(D)$.  So, by Corollary \ref{cor_rankEaD}, there are 
\[\tag{\ref{thm_nrIGS}.1}\label{IGS1}
\Lam!S(r^{n-r},\Lam)
\]
choices for $U_1$.  We multiply this by the number of choices for $U_2$.  
%
By the proof of Theorem \ref{thm_rankEXa}, $\Ub_2$ is a generating set of $\TASA$.  Also, since $|\Ub_2|\leq|U_2|=|U|-|U_1|=\rho_r=\idrank(\TASA)$, it follows that $\Ub_2$ is a minimal idempotent generating set of $\TASA$, and therefore corresponds to a unique graph $\Ga\in\bbT_r$.  We will count the number of ways to choose $U_2$ so that $\Ub_2$ corresponds to $\Ga$.  Consider an edge $(i,j)$ in~$\Ga$.  Then $\ep_{ij}\in\Ub_2$, so there is a unique idempotent $\ve_{ij}\in U_2$ with $\ep_{ij}=\veb_{ij}$.  To specify $\ve_{ij}$, we first choose $\im(\ve_{ij})=\{b_1,\ldots,b_{j-1},b_{j+1},\ldots,b_r\}$, where $b_k\in A_k$ for each $k$.  There are $\lam_1\cdots\lam_{j-1}\lam_{j+1}\cdots\lam_r=\Lam/\lam_j$ choices for $\im(\ve_{ij})$.  Once $\im(\ve_{ij})$ is chosen, $\ve_{ij}$ is restricted by the fact that $a_k\ve_{ij}=b_{ke_{ij}}$ for each $k$.  But the remaining $n-r$ elements of $X\sm A$ may be mapped by $\ve_{ij}$ arbitrarily into the $r-1$ elements of $\im(\ve_{ij})$, and there are $(r-1)^{n-r}$ ways to make these choices.  So the total number of choices for $\ve_{ij}$ is equal to
$
(r-1)^{n-r}\Lam/\lam_j.
$
Since this value depends only on $j$, and since there are $d_\Ga^+(j)$ edges of the form $(i,j)$, taking the product over all edges of $\Ga$ gives a total of
\[\tag{\ref{thm_nrIGS}.2}\label{IGS2}
\prod_{j\in\br} \left[(r-1)^{n-r}\frac{\Lam}{\lam_j}\right]^{d_\Ga^+(j)} 
= 
\big[(r-1)^{n-r}\Lam\big]^{\rho_r}\frac{1}{\lam_1^{d_\Ga^+(1)}\cdots\lam_r^{d_\Ga^+(r)}}
\]
choices for $U_2$ with $\Ub_2$ corresponding to $\Ga$ (noting that $\sum_{j\in\br}d_\Ga^+(j)=\rho_r$).  Summing \eqref{IGS2} over all $\Ga\in\bbT_r$ and multiplying by \eqref{IGS1} gives the result.~\epf

\begin{rem}
Theorem \ref{thm_nrIGS} is also valid if $r=n$, giving $|\bbT_n|$ generating sets for $\E_X=\la E(\T_X)\ra$ of size $1+\rho_n$, in agreement with Theorem \ref{thm_GaU}.  When $r=2$, the given expression reduces to $\Lam^2\Lam!S(2^{n-2},\Lam)$.  
\end{rem}

\section{Ideals of $\RegTXa$}\label{sect:ideals}

In this final section, we consider the ideals of $P=\RegTXa$.  In particular, we show that each of the proper ideals is idempotent generated, and we calculate the rank and idempotent rank, showing that these are equal.  Again, the problems of this section have been solved in the case $r=n$ and are trivial if $r=1$, so we continue to assume that $1<r<n$.  We first state the corresponding result for full transformation semigroups; for convenience, we state it in the context of $\T_A$.

\ms
\begin{thm}[Howie and McFadden \cite{Howie1990}]\label{thm_TX_ideals}
The ideals of $\T_A$ are precisely the sets
\[
I_m=\bigcup_{j\in \bm}D_j=\set{f\in\T_A}{\rank(f)\leq m} \qquad\text{for $1\leq m\leq r$,}
\]
and they form a chain: $I_1\sub\cdots\sub I_r$.  If $m<r$, then $I_m=\la E(D_m)\ra$ is generated by the idempotents in its top $\gD$-class, and
\[\epfreseq
\rank(I_m)=\idrank(I_m)=\begin{cases}
S(r,m) &\text{if $1<m<r$}\\
r &\text{if $m=1$.}
\end{cases}
\]
\end{thm}

The next result is a strengthening Lemma \ref{lemma_r^n-r}.

\newpage
\begin{lemma}\label{lemma_mnrSrm}
If $2\leq m\leq r$, then 
\ms
\bmc2
\itemit{i} $m^{n-r}\geq\Lam_I$ for all $I\in{\br \choose m}$,
\itemit{ii} $\ds{m^{n-r}S(r,m) \geq \sum_{I\in{\br \choose m}}\Lam_I}$.
\emc
\end{lemma}

\pf 
Let $I\in{\br\choose m}$.  Since $\lam_j\geq1$ for all $j\in\br\sm I$, $\sum_{i\in I}\lam_i\leq n-r+m$.  As in the proof of Lemma~\ref{lemma_r^n-r}, we deduce that $\Lam_I\leq(n-r+m)^m/m^m$.  So it suffices to prove that $(n-r+m)^m/m^m\leq m^{n-r}$, which is equivalent to
\[\tag{\ref{lemma_mnrSrm}.1}\label{mnrSrm.1}
(n-r+m)^m\leq m^{n-r+m}.
\]
Note that $n-r+m>m$, so again, as in the proof of Lemma \ref{lemma_r^n-r}, \eqref{mnrSrm.1} is true unless $n-r+m=3$ and $m=2$.  But in this exceptional case, we have $r=n-1$ and $m=2$ so that, without loss of generality, $(\lam_1,\ldots,\lam_r)=(2,1,\ldots,1)$, giving $\Lam_I\leq 2=m^{n-r}$.  This completes the proof of (i).  For (ii), we have
\[
\sum_{I\in{\br\choose m}}\Lam_I \leq {r\choose m}m^{n-r} \leq S(r,m)m^{n-r},
\]
where we have used (i) and the fact that $S(r,m)\geq{r\choose m}$. \epf

\begin{rem}
It follows from Theorem \ref{inflation_thm} and Lemma \ref{lemma_mnrSrm} that each $\gHh$- and $\gDh=\gDa$-class of $P=\RegTXa$ not contained in $D_1^a$ is at least as ``tall'' as it is ``wide''; that is, if $C$ is such a class, then $|C/{\gRa}|\geq|C/{\gLa}|$.
\end{rem}

\ms
\begin{thm}\label{thm_RegTXa_ideals}
The ideals of $P=\RegTXa$ are precisely the sets 
\[
I_m^a=\bigcup_{j\in\bm}D_j^a=\set{f\in P}{\rank(f)\leq m} \qquad\text{for $1\leq m\leq r$,}
\]
and they form a chain: $I_1^a\sub\cdots\sub I_r^a$.  If $m<r$, then $I_m^a = \la \Ea(D_m^a) \raa$ is generated by the idempotents in its top $\gDa$-class, and
\[
\rank(I_m^a)=\idrank(I_m^a)=\begin{cases}
m^{n-r}S(r,m) &\text{if $1<m<r$}\\
n &\text{if $m=1$.}
\end{cases}
\]
\end{thm}

\pf More generally, it may easily be checked that if the $\gJ$-classes of a semigroup $S$ form a chain, $J_1<\cdots<J_q$, then the ideals of $S$ are precisely the sets $I_p=J_1\cup\cdots\cup J_p$ for $1\leq p\leq q$.  Now suppose $m<r$, and let $f\in I_m^a$ be arbitrary.  By Theorem \ref{thm_TX_ideals}, $\fb=h_1\cdots h_k$ for some $h_1,\ldots,h_k\in E(D_m)$.  By Lemma~\ref{lemma_EDm_EDma}, we may choose $e_1,\ldots,e_k\in\Ea(D_m^a)$ such that $\eb_i=h_i$ for each $i$.  Now put $g=e_1\star\cdots\star e_k\in\la\Ea(D_m^a)\raa$.  Then $\gb=\fb$, so by Lemma \ref{lemma_e1e2}, there exist $e_0,e_{k+1}\in\Ea(D_m^a)$ such that $f=e_0\star g\star e_{k+1}\in\la\Ea(D_m^a)\raa$.  

We now prove the statement about rank and idempotent rank.  Note that $I_1^a=D_1^a=D_1$ is an $n$-element right zero semigroup, so the result is trivial for $m=1$; see also Remark \ref{rem_Rank_RegTXA}.  So we assume $1<m<r$ from now on.  More generally, if $J$ is a maximal $\gJ$-class of a finite semigroup $S$, and if the $\gR$-classes contained in~$J$ are $R_1,\ldots,R_q$, then any generating set for $S$ must intersect each $R_p$ non-trivially (for example, this follows from \cite[Exercise~12, p98]{Howie} or from \emph{stability} \cite[Definition~A.2.1]{RSbook}).  In particular, it follows from Theorem \ref{inflation_thm}(vi) that $\rank(I_m^a)\geq m^{n-r}S(r,m)$.  To complete the proof, it suffices to show that there exists $U\sub\Ea(D_m^a)$ with $|U|=m^{n-r}S(r,m)$ and such that $\la U\raa$ contains $D_m^a$.  We now construct such a $U$.

First, let $V\sub E(D_m)\sub\T_A$ be such that $|V|=S(r,m)$ and $I_m=\la V\ra$.  Fix some $v\in V$, write $\im(v)=\{a_{i_1},\ldots,a_{i_m}\}$, and put $I=\{i_1,\ldots,i_m\}$.  Then $H_v\phi^{-1}$ is an $\gHh$-class of $P$, and is an $m^{n-r}\times\Lam_I$ rectangular group.  Put $B_v=\Ea(H_v\phi^{-1})$, so that $B_v$ is a $m^{n-r}\times\Lam_I$ rectangular band.  Since $m^{n-r}\geq\Lam_I$, we see by Lemma \ref{lemma_IJ} that $\rank(B_v)=\idrank(B_v)=m^{n-r}$, so we may choose some $U_v\sub B_v$ with $|U_v|=m^{n-r}$ and $B_v=\la U_v\raa$, noting that $\ub=v$ for all $u\in U_v$.  
Note that for any $h\in D_m^a$ with $\ker(\hb)=\ker(v)$, $U_v$ contains some $u$ with $\ker(u)=\ker(h)$; a similar statement holds for images.
Now put $U=\bigcup_{v\in V}U_v$, noting that $U\sub\Ea(D_m^a)$ and $|U|=m^{n-r}S(r,m)$.  Let $f\in D_m^a$ be arbitrary, and consider an expression $\fb=v_1\cdots v_k$, where $v_1,\ldots,v_k\in V$.  Note that, since $\rank(v_j)=m$ for each $j\in\bk$, we have $\ker(v_1)=\ker(\fb)$ and $\im(v_k)=\im(\fb)$.  For each $j\in\bk$, we choose some $u_j\in U_{v_j}$, but we make these choices so that $\ker(u_1)=\ker(f)$ and $\im(u_k)=\im(f)$.  Put $g=u_1\star\cdots\star u_k$.  Then $\rank(g)=\rank(\gb)=\rank(\fb)=\rank(f)$ so, since $\rank(u_j)=\rank(\ub_j)=m$ for each $j$, we see that $\ker(g)=\ker(u_1)=\ker(f)$ and $\im(g)=\im(u_k)=\im(f)$.    Together with $\fb=\gb$, this shows that $f=g\in\la U\raa$, and completes the proof. \epf

\begin{rem}\label{rem_Rank_RegTXA2}
Again, we note the similarity between Theorem \ref{thm_RegTXa_ideals} and \cite[Theorem 4.4]{SS2013}, where it is shown that the proper ideals, there denoted $Q(F;m)$, of $\RegTXA$ are idempotent generated, and have rank and idempotent rank equal to $m^{n-r}S(r,m)$.  See also Remark \ref{rem_Rank_RegTXA}.  We also note that an alternative approach exists to tackle problems such as those we addressed in this section; namely, making use of the general results of Gray \cite{Gray2007,Gray2008} on (idempotent) rank in completely $0$-simple semigroups.
\end{rem}

\subsection*{Acknowledgement}

We kindly thank the referee for his/her careful reading of the manuscript and for some helpful suggestions.

\footnotesize
\def\bibspacing{-1.1pt}
\bibliography{variants_biblio}
\bibliographystyle{plain}
\end{document}